\documentclass{ws-ijbc}
\usepackage{ws-rotating}     
\usepackage{graphicx}
\usepackage{epstopdf}

\usepackage{amsfonts, amsbsy, amssymb, amsmath, commath, graphicx, float, footnote, natbib, subfigure, color, url, 
titlesec}
\usepackage{comment}
\usepackage[english]{babel}

\begin{document}

\bibliographystyle{plain}

\markboth{Balibrea-Iniesta et al.}{}

\title{\Large{\textbf{LAGRANGIAN DESCRIPTORS FOR STOCHASTIC DIFFERENTIAL EQUATIONS: A TOOL FOR REVEALING THE PHASE PORTRAIT OF STOCHASTIC DYNAMICAL SYSTEMS}}}

\author{FRANCISCO BALIBREA-INIESTA\footnote{Instituto de Ciencias Matem\'aticas, CSIC-UAM-UC3M-UCM, Spain, \textit{francisco.balibrea@icmat.es}.}}

\address{Instituto de Ciencias Matem\'aticas, CSIC-UAM-UC3M-UCM,\\ C/ Nicol\'as Cabrera 15, Campus de Cantoblanco UAM\\
Madrid, 28049, Spain\\
francisco.balibrea@icmat.es}

\author{CARLOS LOPESINO\footnote{Instituto de Ciencias Matem\'aticas, CSIC-UAM-UC3M-UCM, Spain, \textit{carlos.lopesino@icmat.es}.}}

\address{Instituto de Ciencias Matem\'aticas, CSIC-UAM-UC3M-UCM,\\ C/ Nicol\'as Cabrera 15, Campus de Cantoblanco UAM\\
Madrid, 28049, Spain\\
carlos.lopesino@icmat.es}

\author{STEPHEN WIGGINS\footnote{School of Mathematics, University of Bristol, United Kingdom, \textit{s.wiggins@bristol.ac.uk}.}}

\address{School of Mathematics, University of Bristol,\\ Bristol BS8 1TW, United Kingdom \\
s.wiggins@bristol.ac.uk}

\author{ANA M. MANCHO\footnote{Instituto de Ciencias Matem\'aticas, CSIC-UAM-UC3M-UCM, Spain, \textit{a.m.mancho@icmat.es}.}}

\address{Instituto de Ciencias Matem\'aticas, CSIC-UAM-UC3M-UCM,\\ C/ Nicol\'as Cabrera 15, Campus de Cantoblanco UAM\\
Madrid, 28049, Spain\\
a.m.mancho@icmat.es}

\maketitle

\begin{abstract}
\noindent
In this paper we  introduce a new technique for depicting the phase portrait of stochastic differential equations. Following previous work for deterministic systems, we represent the phase space by means of a generalization of the method of Lagrangian descriptors to stochastic differential equations. Analogously to the deterministic differential equations setting, the Lagrangian descriptors graphically provide the distinguished trajectories and hyperbolic structures arising within the stochastic dynamics, such as  random fixed points and their stable and unstable manifolds. We analyze the sense in which structures form barriers to transport in stochastic systems. We apply the method to several benchmark examples where the deterministic phase space structures are well-understood. In particular, we apply our method to the noisy saddle, the stochastically forced Duffing equation, and the stochastic double gyre model that is a benchmark for analyzing fluid transport.
\newline
\\
\textit{Keywords:} Stochastic differential equation, Lagrangian descriptor, distinguished trajectories, hyperbolicity, random fixed point, stable and unstable manifolds.
\end{abstract}

\section{Introduction}
\label{sec:intro}

In recent years Lagrangian descriptors (LD) have been shown to be a useful technique for discovering  phase space structure in both autonomous and nonautonomous dynamical systems. In \cite{chaos} Lagrangian descriptors (also referred to in the literature as the ``M function'') were used to discover hyperbolic trajectories and their stable and unstable manifolds in aperiodically time dependent vector fields. In particular, this paper considers ``important hyperbolic trajectories'', referred to as
{\em distinguished trajectories}, which build on the well-known idea of distinguished {\em hyperbolic} trajectory,  discussed earlier in \cite{kayo}.

Lagrangian descriptors can easily be applied to the analysis of velocity fields defined as data sets. Early work demonstrating this was concerned with transport associated with the Kuroshio current, described in \cite{prl,jfm}. This work set the stage for further geophysical transport studies described in, for example, \cite{amism11, alvaro2}. An application concerned with determining the connection between
coherent structures and the saturation of a nonlinear dynamo is described  in \cite{rempel}.  More recently, Lagrangian descriptors have been used to analyze issue related to the search strategy for the Malaysian airliner MH370 \cite{GMWM15} as well as to an understanding of events surrounding the Deepwater Horizon oil spill
(\cite{mmw14}) and, more recently, to an analysis of the real time development
of an oil spill in the Canary Islands \cite{GRMCW15}.

In \cite{cnsns} a general assessment of Lagrangian descriptors was carried out.
Different Lagrangian descriptors were considered and applied to a variety of benchmark examples containing known phase space structures. While most of the previous applications of Lagrangian descriptors had been to two dimensional, time dependent flows, it was also shown that dimensionality is not a restriction to the application of the method. In particular, transport associated with the three dimensional time dependent Hill's spherical vortex was considered. The computational requirements of Lagrangian descriptors were also considered and compared with those of finite time Lyapunov exponents (FTLEs).

The applications of Lagrangian descriptors mentioned above are in areas of fluid mechanics. However, the method of Lagrangian descriptors applies to  general vector fields in $n$ dimensions, and there is no obstacle to using Lagrangian descriptors for applications in a higher dimensional setting. Recently, this has been illustrated by a series of applications of Lagrangian descriptors to problems in chemistry by Hernandez and co-workers. In particular, they have applied Lagrangian descriptors to a study of chemical reactions under external time-dependent driving in \cite{CH15}, barrierless reactions in \cite{JH16, JH2016}, and ketene isomerization in \cite{CH16}. Finally, we remark that the use of  Lagrangian descriptors for visualizing phase space structures in complex dynamical systems is described in \cite{MCWGM15}.

In this paper we continue the development of the method of Lagrangian descriptors by carrying it out in the setting of stochastic differential equations. This provides us with a computational tool for revealing the phase space structure of stochastic dynamical systems. We will discuss in detail what this means by considering explicit examples and comparing them with their deterministic counterparts.

This paper is organized as follows: in Section \ref{sec:pc} we give a brief summary of preliminary concepts as well as 
set the general framework of stochastic differential equations. In Section \ref{sec:SLD} we define Stochastic Lagrangian 
descriptor (SLD) as a method to reveal the phase space of a stochastic system and in Section \ref{sec:num} we describe  the numerical method used to implement SLD. Finally in Section \ref{sec:examp}, we compute the SLD by 
using some benchmark examples: the noisy saddle system, the stochastic versions of the Duffing equation, and the stochastic
Double Gyre.

\section{Preliminary concepts}
\label{sec:pc}

The general framework in which  Lagrangian descriptors (in their stochastic version) are defined, requires us to carefully describe the nature of solutions of a stochastic differential equation (SDE). Therefore, to begin we consider a general system of SDEs expressed in  differential form as follows:

\begin{equation}
\label{SDE}
dX_{t} = b(X_{t},t)dt + \sigma (X_{t},t)dW_{t}, \quad t \in \mathbb{R},
\end{equation}

\noindent
where $b(\cdot) \in C^{1}(\mathbb{R}^{n}\times \mathbb{R})$ is the deterministic part, $\sigma (\cdot) \in C^{1}(\mathbb{R}^{n}\times \mathbb{R})$ is the random forcing, $W_{t}$ is a Wiener process (also called Brownian motion) whose definition is given later, and $X_{t}$ is the solution of the equation. All these functions take values in $\mathbb{R}^{n}$.

As the notion of solution of a SDE is closely related with the Wiener process, we  state what is meant by $W(\cdot )$. This definition is given in \cite{duan15}, and this reference serves to provide the background for all of the notions in this section.
Also, throughout this article we will use $\Omega$ to denote the probability space where the Wiener process is defined.

\begin{definition}
\label{Wiener}
A real valued stochastic Wiener or Brownian process $W(\cdot)$ is a stochastic process defined in a probability space $(\Omega , {\cal F},{\cal P})$ which satisfies
\begin{enumerate}
  \item[(i)] $W_0 = 0$ (standard Brownian motion),
  \item[(ii)] $W_t - W_s$  follows a Normal distribution $N(0,t-s)$ for all $t\geq s \geq 0$,
  \item[(iii)] for all time $0 < t_1 < t_2 < ... < t_n$, the random variables $W_{t_1}, W_{t_2} - W_{t_1},... , W_{t_n} - W_{t_{n-1}}$ are independent (independent increments).
\end{enumerate}
Moreover, $W(\cdot)$ is a real valued two-sided Wiener process if conditions (ii) and (iii) change into
\begin{enumerate}
  \item[(ii)] $W_t - W_s$ follows a Normal distribution $N(0,|t-s|)$ for all $t, s \in \mathbb{R}$,
  \item[(iii)] for all time $t_1 , t_2 , ... , t_{2n} \in \mathbb{R}$ such that the intervals $\lbrace (t_{2i-1},t_{2i}) \rbrace_{i=1}^{n}$ are non-intersecting between them\footnote{With the notation $(t_{2i-1},t_{2i})$ we refer to the interval of points between the values $t_{2i-1}$ and $t_{2i}$, regardless the order of the two extreme values. Also with the asertion we impose that every pair of intervals of the family $\lbrace (t_{2i-1},t_{2i}) \rbrace_{i=1}^{n}$ has an empty intersection, or alternatively that the union $\bigcup_{i=1}^{n}(t_{2i-1},t_{2i})$ is conformed by $n$ distinct intervals over $\mathbb{R}$.}, the random variables $W_{t_1}-W_{t_2}, W_{t_3} - W_{t_4},... , W_{t_{2n-1}} - W_{t_{2n}}$ are independent.
\end{enumerate}
\end{definition}

As it is mentioned in the introduction, the purpose of this article is to develop the method of Lagrangian descriptors for SDEs. This method of Lagrangian descriptors  has been developed for deterministic differential equations whose temporal domain is $\mathbb{R}$. In this sense it is natural to work with two-sided solutions as well as two-sided Wiener processes. Henceforth, every Wiener process $W(\cdot )$ considered in the this article will be of this form.

Given that any Wiener process $W(\cdot )$ is a stochastic process, by definition this is a family of random real variables $\lbrace W_{t}, t\in \mathbb{R} \rbrace$ in such a way that for each $\omega \in \Omega$ there exists a mapping
$$ t \longmapsto W_{t}(\omega )$$
known as the trajectory of a Wiener process.

Analogously to the Wiener process, the solution $X_{t}$ of the SDE \eqref{SDE} is also a stochastic process. In particular, it is a family of random variables $\lbrace X_{t}, t\in \mathbb{R} \rbrace$ such that for each $\omega \in \Omega$, the trajectory of $X_{t}$ satisfies

\begin{equation}
\label{Xt}
t \longmapsto X_{t}(\omega ) = X_{0}(\omega ) + \int_{0}^{t} b(X_{s}(\omega ), s)ds + \int_{0}^{t} \sigma (X_{s}(\omega ), s)dW_{s}(\omega ),
\end{equation}

\noindent
where $X_{0}:\Omega \rightarrow \mathbb{R}^{n}$ is the initial condition. In addition, as $b(\cdot)$ and $\sigma(\cdot)$ are smooth functions, they are locally Lipschitz and this leads to existence and pathway uniqueness of a local, continuous solution (see \cite{duan15}). That is if any two stochastic processes $X^1$ and $X^2$ are local solutions in time of SDE \eqref{SDE}, then $X^1_t(\omega) = X^2_t(\omega)$ over a time interval $t \in (t_{i},t_{f})$ and for almost every $\omega \in \Omega$.

At each instant of time $t$, the deterministic integral $\int_{0}^{t} b(X_{s}(\omega ))ds$ is defined by the usual Riemann integration scheme since $b$ is assumed to be a differentiable function. However,  the stochastic integral term is chosen to be defined by the It\^{o} integral scheme:

\begin{equation}
\label{Ito}
\int_{0}^{t} \sigma (X_{s}(\omega ),s)dW_{s}(\omega ) = \lim_{N \rightarrow \infty} \sum_{i=0}^{N-1} \sigma (X_{i\frac{t}{N}}(\omega ), it/N) \cdot \left[ W_{(i+1)\frac{t}{N}}(\omega ) - W_{i\frac{t}{N}}(\omega ) \right].
\end{equation}

\noindent
This scheme will also facilitate the implementation of a numerical method for computing approximations for the solution $X_{t}$ in the next section.

Once the notion of solution, $X_{t}$, of a SDE (\ref{SDE}) is established, it is natural to ask if the same notions and ideas familar from the study of deterministic differential equations  from the dynamical systems point of view are still valid for SDEs. In particular, we want to consider the notion of hyperbolic trajectory and its stable and unstable manifolds in the context of SDEs. We  also want to consider how such notions would manifest themselves in the context of {\em phase space transport} for SDEs, and the stochastic Lagrangian descriptor will play a key role in considering these questions from a practical point of view.

We first discuss the notion of an invariant set for a SDE.
In the deterministic case the simplest possible invariant set is a single trajectory of the differential equation. More precisely, it is the set of points through which a solution passes.  Building on this construction, an invariant set is a collection of trajectories of different solutions. This is the most basic way to characterize the invariant sets with respect to a determinsitic differential equation of the form

\begin{equation}
\label{deterministic_system}
\dot{x} = f(x,t), \quad x \in \mathbb{R}^{n}, \quad t \in \mathbb{R}.
\end{equation}

\noindent
For verifying  the invariance of such sets the solution mapping generated by the vector field is used.  For  deterministic autonomous systems  these are referred to as {\em flows} (or `'dynamical systems'') and for deterministic nonautonomous systems they are referred to as {\em processes}. The formal definitions are in Appendix 1 and can also be found in \cite{kloe11}. 

A similar notion of solution mapping for SDEs is introduced using the notion of a random dynamical system $\varphi$ (henceforth referred to as RDS) in the context of SDEs. This function $\varphi$ is also a solution mapping of a SDE that satisfies several conditions, but compared with the solution mappings in the deterministic case, this RDS depends on an extra argument which is the random variable $\omega \in \Omega$. Furthermore the random variable $\omega$ evolves with respect to $t$ by means of a dynamical system $\lbrace \theta_{t} \rbrace_{t \in \mathbb{R}}$ defined over the probability space $\Omega$.

\begin{definition}
\label{rds}
Let $\lbrace \theta_{t} \rbrace_{t \in \mathbb{R}}$ be a measure-preserving\footnote{Given the probability measure $\mathcal{P}$ associated with the space $(\Omega , \mathcal{F},\mathcal{P})$, this remains invariant under the dynamical system $\lbrace \theta_{t} \rbrace_{t \in \mathbb{R}}$. Formally, $\theta_{t}\mathcal{P} = \mathcal{P}$ for every $t \in \mathbb{R}$. This statement means that $\mathcal{P}(B)=\mathcal{P}(\theta_{t}(B))$ for every $t \in \mathbb{R}$ and every subset $B \in \mathcal{F}$. Indeed for any dynamical system $\lbrace \theta_{t} \rbrace_{t \in \mathbb{R}}$ defined over the same probability space $\Omega$ as a Wiener process $W(\cdot )$, we have the equality $W_{s}(\theta_{t}\omega ) = W_{t+s}(\omega )-W_{t}(\omega )$ which implies that $dW_{s}(\theta_{t}\omega ) = dW_{t+s}(\omega )$ for every $s,t \in \mathbb{R}$ (see \cite{duan15} for a detailed explanation).} dynamical system defined over $\Omega$, and let $\varphi : \mathbb{R} \times \Omega \times \mathbb{R}^{n} \rightarrow \mathbb{R}^{n}$ be a measurable mapping such that $(t,\cdot , x) \mapsto \varphi (t,\omega ,x)$ is continuous for all $\omega \in \Omega$, and the family of functions $\lbrace \varphi (t,\omega ,\cdot ): \mathbb{R}^{n} \rightarrow \mathbb{R}^{n} \rbrace$ has the cocycle property:
$$ \varphi (0,\omega ,x)=x \quad \text{and} \quad \varphi (t+s,\omega ,x) = \varphi(t,\theta_{s}\omega,\varphi (s,\omega ,x)) \quad \text{for all } t,s \in \mathbb{R}, \text{ } x \in \mathbb{R}^{n} \text{ and } \omega \in \Omega .$$
Then the mapping $\varphi$ is a random dynamical system with respect to the stochastic differential equation
$$dX_{t} = b(X_{t})dt + \sigma (X_{t})dW_{t}$$
if $\varphi (t,\omega ,x)$ is a solution of the equation.
\end{definition}

Analogous to the deterministic case, the definition of invariance with respect to a SDE can be characterized in terms of a RDS. This is an important topic in our consideration of stochastic Lagrangian descriptors. Now we introduce an example of a SDE for which the analytical expression of the RDS is obtained. This will be a benchmark example in our development of stochastic Lagrangian descriptors their relation to stochastic invariant manifolds.

\begin{example}{(Noisy saddle point)}

For the stochastic differential equation

\begin{equation}
\label{noisy_saddle}
\begin{cases} dX_{t} = X_{t}dt + dW_{t}^{1} \\ dY_{t} = -Y_{t}dt + dW_{t}^{2} \end{cases}
\end{equation}

\noindent
where $W_{t}^{1}$ and $W_{t}^{2}$ are two different Wiener processes, the solutions take the expressions

\begin{equation}
\label{noisy_saddle_solutions}
X_{t} = e^{t} \left( X_{0}(\omega ) + \int_{0}^{t}e^{-s}dW_{s}^{1}(\omega ) \right) \quad , \quad Y_{t} = e^{-t} \left( Y_{0}(\omega ) + \int_{0}^{t}e^{s}dW_{s}^{2}(\omega ) \right)
\end{equation}

\noindent
and therefore the random dynamical system $\varphi$ takes the form

\begin{equation}
\label{noisy_saddle_RDS}
\begin{array}{ccccccccc} \varphi : & & \mathbb{R} \times \Omega \times \mathbb{R}^{2} & & \longrightarrow & & \mathbb{R}^{2} & & \\ & & (t,\omega ,(x,y)) & & \longmapsto & & \left( \varphi_{1}(t,\omega ,x),\varphi_{2}(t,\omega ,y) \right) & = & \left( e^{t} \left( x + \int_{0}^{t}e^{-s}dW_{s}^{1}(\omega ) \right) , e^{-t} \left( y + \int_{0}^{t}e^{s}dW_{s}^{2}(\omega ) \right) \right) . \end{array}
\end{equation}
\end{example}

Notice that this last definition (\ref{rds}) is expressed in terms of SDEs with time-independent coefficients $b,\sigma$. For more general SDEs we introduce a definition of nonautonomous RDS in Appendix 2, inspired in Definition \ref{process} of processes for deterministic nonautonomous differential equations in Appendix 1. However, for the remaining examples considered in this article we make use of the already given definition of RDS.

Once we have the notion of RDS, it can be used to describe and detect geometrical structures and determine their influence on the  dynamics of trajectories. Specifically, in clear analogy with the deterministic case, we focus on those trajectories whose expressions do not depend explicitly on time $t$, which are referred as {\em random fixed points}. Moreover, their stable and unstable manifolds, which may also depend on the random variable $\omega$, are also objects of interest due to their influence on the dynamical behavior of nearby trajectories. Both types of objects are invariant. Therefore we describe a characterization of invariant sets with respect to a SDE by means of an associated RDS.

\begin{definition}
\label{invariant_set}
A non empty collection $M : \Omega \rightarrow \mathcal{P}(\mathbb{R}^{n})$, where $M(\omega ) \subseteq \mathbb{R}^{n}$ is a closed subset for every $\omega \in \Omega$, is called an invariant set for a random dynamical system $\varphi$ if

\begin{equation}
\label{invariance}
\varphi (t,\omega ,M(\omega )) = M(\theta_{t}\omega ) \quad \text{for every } t \in \mathbb{R} \text{ and every } \omega \in \Omega.
\end{equation}

\end{definition}

Again, we return to the noisy saddle (\ref{noisy_saddle}), which  is an example of a SDE for which several invariant sets can be easily characterized by means of its corresponding RDS.

\begin{example}{\textbf{(Noisy saddle point)}}
For the stochastic differential equations

\begin{equation}
\begin{cases} dX_{t} = X_{t}dt + dW_{t}^{1} \\ dY_{t} = -Y_{t}dt + dW_{t}^{2} \end{cases}
\end{equation}

\noindent
where $W_{t}^{1}$ and $W_{t}^{2}$ are two different Wiener processes, the solution mapping $\varphi$ is given by  the following  expression

$$\begin{array}{ccccccccc} \varphi : & & \mathbb{R} \times \Omega \times \mathbb{R}^{2} & & \longrightarrow & & \mathbb{R}^{2} & & \\ & & (t,\omega ,(x,y)) & & \longmapsto & & (\varphi_{1}(t,\omega ,x),\varphi_{2}(t,\omega ,y)) & = & \left( e^{t} \left( x + \int_{0}^{t}e^{-s}dW_{s}^{1}(\omega ) \right) , e^{-t} \left( y + \int_{0}^{t}e^{s}dW_{s}^{2}(\omega ) \right) \right) . \end{array}$$

\noindent
Notice that this is a decoupled random dynamical system. There exists a solution whose components do not depend on variable $t$ and are convergent for almost every $\omega \in \Omega$ as a consequence of the properties of Wiener processes (see \cite{duan15}). This solution has the form:

$$\tilde{X}(\omega) = (\tilde{x}(\omega ),\tilde{y}(\omega )) = \left( - \int_{0}^{\infty}e^{-s}dW_{s}^{1}(\omega ) , \int_{-\infty}^{0}e^{s}dW_{s}^{2}(\omega ) \right) .$$

\noindent
Actually, $\tilde{X}(\omega )$ is a solution because it satisfies the invariance property that we now verify:

\begin{equation}
\label{invariance_x}
\begin{array}{ccl}
\varphi_{1} (t,\omega ,\tilde{x}(\omega )) & = & \displaystyle{e^{t}\left( -\int_{0}^{+\infty}e^{-s}dW^{1}_{s}(\omega ) 
+ \int_{0}^{t}e^{-s}dW^{1}_{s}(\omega ) \right) } = \displaystyle{-\int_{t}^{+\infty}e^{-(s-t)}dW^{1}_{s}(\omega ) }\\ 
& = & \displaystyle{-\int_{0}^{+\infty}e^{-t'}dW^{1}_{t'+t}(\omega) = - 
\int_{0}^{+\infty}e^{-t'}dW^{1}_{t'}(\theta_{t}\omega ) = \tilde{x}(\theta_{t}\omega )} \quad \text{by means of } 
t'=s-t,
\end{array}
\end{equation}

\begin{equation}
\label{invariance_y}
\begin{array}{lll}
\varphi_{2} (t,\omega ,\tilde{y}(\omega )) & = & \displaystyle{e^{-t}\left( \int_{-\infty}^{0}e^{s}dW^{2}_{s}(\omega ) + 
\int_{0}^{t}e^{s}dW^{2}_{s}(\omega ) \right)  = \int_{-\infty}^{t}e^{s-t}dW^{2}_{s}(\omega ) = 
\int_{-\infty}^{0}e^{t'}dW^{2}_{t'+t}(\omega )} \\  & = & \displaystyle{ 
\int_{-\infty}^{0}e^{t'}dW^{2}_{t'}(\theta_{t}\omega ) = \tilde{y}(\theta_{t}\omega )} \quad \text{by means of } t'=s-t.
\end{array}
\end{equation}

\noindent
This implies that $\varphi (t,\omega ,\tilde{X}(\omega )) = \tilde{X}(\theta_{t} \omega)$ for every $t \in \mathbb{R}$ and every $\omega \in \Omega$. Therefore $\tilde{X}(\omega )$ satisfies the invariance property (\ref{invariance}). This conclusion comes from the fact that $\tilde{x}(\omega )$ and $\tilde{y}(\omega )$ are also invariant under the components $\varphi_{1}$ and $\varphi_{2}$, in case these are seen as separate RDSs defined over $\mathbb{R}$ (see (\ref{invariance_x}) and (\ref{invariance_y}), respectively).

Due to its independence with respect to the time variable $t$, it is said that $\tilde{X}(\omega )$ is a random fixed point of the SDE (\ref{noisy_saddle}), or more commonly a stationary orbit. As the trajectory of $\tilde{X}(\omega )$ (and separately its components $\tilde{x}(\omega )$ and $\tilde{y}(\omega )$) is proved to be an invariant set, it is straightforward to check that the two following subsets of $\mathbb{R}^{2}$,

$$\mathcal{S}(\omega ) = \lbrace (x,y) \in \mathbb{R}^{2} : x = \tilde{x}(\omega ) \rbrace \quad , \quad \mathcal{U}(\omega ) = \lbrace (x,y) \in \mathbb{R}^{2} : y = \tilde{y}(\omega ) \rbrace $$

\noindent
are also invariant with respect to the RDS $\varphi$. Similarly to the deterministic setting, these are referred to as the stable and unstable manifolds of the stationary orbit respectively. Additionally, in order to prove the separating nature of these two manifolds and the stationary orbit with respect to their nearby trajectories, let's consider any other solution $(\overline{x}_{t},\overline{y}_{t})$ of the noisy saddle with initial conditions at time $t=0$,

$$\overline{x}_{0} = \tilde{x}(\omega ) + \epsilon_{1}(\omega ) , \quad \overline{y}_{0} = \tilde{y}(\omega ) + \epsilon_{2}(\omega ), \quad \text{being } \epsilon_{1}(\omega ), \epsilon_{2}(\omega ) \text{ two random variables.}$$

\noindent
If the corresponding RDS $\varphi$ is applied to compare the evolution of this solution $(\overline{x}_{t},\overline{y}_{t})$ and the stationary orbit, there arises an exponential dichotomy:

$$ (\overline{x}_{t},\overline{y}_{t}) - (\tilde{x}(\theta_{t}\omega ),\tilde{y}(\theta_{t}\omega )) = \varphi (t,\omega ,(\overline{x}_{0},\overline{y}_{0})) - \varphi (t,\omega ,(\tilde{x}(\omega ),\tilde{y}(\omega ))) $$

$$= \left( e^{t}\left[ \overline{x}_{0} + \int_{0}^{t}e^{-s}dW_{s}^{1}(\omega ) - \tilde{x}(\omega ) - \int_{0}^{t}e^{-s}dW_{s}^{1}(\omega ) \right] , e^{-t}\left[ \overline{y}_{0} + \int_{0}^{t}e^{s}dW_{s}^{2}(\omega ) - \tilde{y}(\omega ) - \int_{0}^{t}e^{s}dW_{s}^{2}(\omega ) \right] \right) $$

\begin{equation}
\label{dichotomy}
= \left( e^{t} \left( \tilde{x}(\omega )+\epsilon_{1}(\omega )-\tilde{x}(\omega ) \right) ,e^{-t} \left( \tilde{y}(\omega )+\epsilon_{2}(\omega )-\tilde{y}(\omega ) \right) \right) = \left( \epsilon_{1}(\omega )e^{t},\epsilon_{2}(\omega )e^{-t} \right) .
\end{equation}

\noindent
Considering that $(\overline{x}_{t},\overline{y}_{t})$ is different from $(\tilde{x}(\omega ),\tilde{y}(\omega ))$ then one of the two cases $\epsilon_{1} \not \equiv 0$ or $\epsilon_{2} \not \equiv 0$ holds, let say $\epsilon_{1} \not = 0$ or $\epsilon_{2} \not = 0$ for almost every $\omega \in \Omega$. In the first case, the distance between both trajectories $(\overline{x}_{t},\overline{y}_{t})$ and $(\tilde{x},\tilde{y})$ increases at an exponential rate in positive time:

\begin{equation}
\label{eq_1}
\norm{\varphi (t,\omega ,(\overline{x}_{t},\overline{y}_{t}))-\varphi (t,\omega ,(\tilde{x},\tilde{y}))} \geq |\epsilon_{1}(\omega )e^{t}| \longrightarrow + \infty \quad \text{when }  t \rightarrow + \infty \text{ and for a.e. } \omega \in \Omega .
\end{equation}

\noindent
Similarly to this case, when the second option holds the distance between both trajectories increases at a exponential rate in negative time. It does not matter how close the initial condition $(\overline{x}_{0},\overline{x}_{0})$ is from $(\tilde{x}(\omega ),\tilde{y}(\omega ))$ at the initial time $t=0$. Actually this same exponentially growing separation can be achieved for any other initial time $t \not = 0$. Following these arguments, one can check that the two manifolds $\mathcal{S}(\omega )$ and $\mathcal{U}(\omega )$ also exhibit this same separating behaviour as the stationary orbit. Moreover, we remark  that almost surely the stationary orbit is the only solution whose components are bounded.

\end{example}

These facts highlight the distinguished nature of the stationary orbit (and its manifolds) in the sense that it is an isolated solution from the others. Apart from the fact that $(\tilde{x},\tilde{y})$ ``moves" in a bounded domain for every $t \in \mathbb{R}$, any other trajectory eventually passing through an arbitrary neighborhood of $(\tilde{x},\tilde{y})$ at any given instant of time $t$,  leaves the neighborhood and then separates from the stationary orbit in either positive or negative time. Specifically, this separation rate is exponential for the noisy saddle, just in the same way as for the deterministic saddle.

These features are also observed for the trajectories within the stable and unstable manifolds of the stationary orbit, but in a more restrictive manner than $(\tilde{x},\tilde{y})$. Taking for instance an arbitrary trajectory $(x^{s},y^{s})$ located at $\mathcal{S}(\omega )$ for every $t \in \mathbb{R}$, its first component $x^{s}_{t}=\tilde{x}(\omega )$ remains bounded for almost every $\omega \in \Omega$. In contrast, any other solution passing arbitrarily closed to $(x^{s},y^{s})$ neither being part of $\mathcal{S}(\omega )$ nor being the stationary orbit, satisfies the previous inequality (\ref{eq_1}) and therefore separates from $\mathcal{S}(\omega )$ at an exponential rate for increasing time. With this framework we can now introduce the formal definitions of stationary orbit and invariant manifold.

\begin{definition}
\label{stationary_orbit}
A random variable $\tilde{X}: \Omega \rightarrow \mathbb{R}^{n}$ is called a stationary orbit (random fixed point) for a random dynamical system $\varphi$ if

$$\varphi(t, \omega, \tilde{X}(\omega)) = \tilde{X}(\theta_t\omega), \quad \text{for every } t \in \mathbb{R} \text{ and every } \omega \in \Omega .$$

\end{definition}

Obviously every stationary orbit $\tilde{X}(\omega )$ is an invariant set with respect to a RDS as it satisfies Definition \ref{invariant_set}. Among several definitions of invariant manifolds given in the bibliography (for example \cite{arno98}, \cite{boxl89}, \cite{duan15}), which have different formalisms but share the same philosophy, we choose the one given in \cite{duan15} because it adapts to our example in a very direct way.

\begin{definition}
A random invariant set $M: \Omega \rightarrow \mathcal{P}(\mathbb{R}^{n})$ for a random dynamical system $\varphi$ is called a $C^{k}$-Lipschitz invariant manifold if it can be represented by a graph of a $C^{k}$ Lipschitz mapping ($k \geq 1$)
$$\gamma (\omega , \cdot ):H^{+} \to H^{-}, \quad \text{with direct sum decomposition } H^{+} \oplus H^{-} = \mathbb{R}^{n}$$
$$\text{such that} \quad M(\omega ) = \lbrace x^{+} \oplus \gamma(\omega ,x^{+}) : x^{+} \in H^{+} \rbrace \quad \text{for every } \omega \in \Omega.$$
\end{definition}

This is a very limitated notion of invariant manifold as its formal definition requires the set to be represented by a Lipschitz graph. Anyway, it is consistent with the already established manifolds $\mathcal{S}(\omega )$ and $\mathcal{U}(\omega )$ as these can be represented as the graphs of two functions $\gamma_{s}$ and $\gamma_{u}$ respectively,

\begin{equation}
\begin{array}{ccccc} \gamma_{s} (\omega , \cdot ) & : & span \lbrace \left( \begin{array}{c} 0 \\ 1 \end{array} \right) \rbrace & \longrightarrow & span \lbrace \left( \begin{array}{c} 1 \\ 0 \end{array} \right) \rbrace \\
& & \left( \begin{array}{c} 0 \\ t \end{array} \right) & \longmapsto & \left( \begin{array}{c} \tilde{x}(\omega ) \\ 0 \end{array} \right) \end{array}
\quad \quad \text{and} \quad \quad
\begin{array}{ccccc} \gamma_{u} (\omega , \cdot ) & : & span \lbrace \left( \begin{array}{c} 1 \\ 0 \end{array} \right) \rbrace & \longrightarrow & span \lbrace \left( \begin{array}{c} 0 \\ 1 \end{array} \right) \rbrace \\
& & \left( \begin{array}{c} t \\ 0 \end{array} \right) & \longmapsto & \left( \begin{array}{c} 0 \\ \tilde{y}(\omega ) \end{array} \right) \end{array} \quad .
\end{equation}
\\

\noindent
Actually the domains of such functions $\gamma_{s}$ and $\gamma_{u}$ are the linear subspaces $E^{s}(\omega )$ and $E^{u}(\omega )$, known as the stable and unstable subspaces of the random dynamical system $\Phi (t,\omega )$. This last mapping is obtained from linearizing the original RDS $\varphi$ over the stationary orbit $(\tilde{x},\tilde{y})$. This result serves as an argument to denote $\mathcal{S}(\omega )$ and $\mathcal{U}(\omega )$ as the stable and unstable manifolds of the stationary orbit, not only because these two subsets are invariant under $\varphi$, as one can deduce from (\ref{invariance_x}) and (\ref{invariance_y}), but also due to the dynamical behaviour of their trajectories in a neighborhood of the stationary orbit $\tilde{X}(\omega )$. Hence the important  characteristic  of $\tilde{X}(\omega )=(\tilde{x},\tilde{y})$ is not only its independence with respect to the time variable $t$; but also the fact  that it exhibits  hyperbolic behaviour with respect to its neighboring trajectories. Considering the hyperbolicity of a given solution, as well as in the deterministic context, means considering  the hyperbolicity of the RDS $\varphi$ linearized over such solution. Specifically, the Oseledets' multiplicative ergodic theorem for random dynamical systems (\cite{arno98} and \cite{duan15}) ensures the existence of a Lyapunov spectrum which is necessary to determine whether the stationary orbit $\tilde{X}(\omega )$ is hyperbolic or not. All these issues are well reported in \cite{duan15} and summarized in Appendix 3, including the proof that the noisy saddle (\ref{noisy_saddle}) satisfies the Oseledets' multiplicative ergodic theorem conditions.

Before  implementing the numerical method of Lagrangian descriptors to several examples of SDEs, it is important to remark why random fixed points and their respective stable and unstable manifolds govern the nearby trajectories, and furthermore, how they may influence the dynamics throughout the rest of the domain. These are essential issues in order to describe the global phase space  motion of  solutions of SDEs. However, these questions do not have a simple answer. For instance, in the noisy saddle example (\ref{noisy_saddle}) the geometrical structures arising from the dynamics generated around the stationary orbit are quite similar to the dynamics corresponding to the deterministic saddle point $\lbrace \dot{x}=x,\dot{y}=-y \rbrace$. Significantly, the manifolds $\mathcal{S}(\omega )$ and $\mathcal{U}(\omega )$ of the noisy saddle form two dynamical barriers for other trajectories in the same way that the manifolds $\lbrace x = 0 \rbrace$ and $\lbrace y = 0 \rbrace$ of the deterministic saddle work. This means that for any particular experiment, i.e., for any given $\omega \in \Omega$, the manifolds $\mathcal{S}(\omega )$ and $\mathcal{U}(\omega )$ are determined and cannot be ``crossed'' by other trajectories due to the uniqueness of solutions (remember that the manifolds are invariant under the RDS (\ref{noisy_saddle_RDS}) and are comprised of an infinite family of solutions). Also by considering the exponential separation rates reported in (\ref{eq_1}) with the rest of trajectories, the manifolds $\mathcal{S}(\omega )$ and $\mathcal{U}(\omega )$ divide the plane $\mathbb{R}^{2}$ of initial conditions into four qualitatively distinct dynamical regions; therefore providing a phase portrait representation. 

Nevertheless it remains to show that such analogy can be found between other SDEs and their corresponding non-noisy deterministic differential equations\footnote{Otherwise if nonlinearity is dominating the behavior of the terms in equation (\ref{SDE}) then the correspondence between the manifolds for $\Phi (t, \omega )$ to the manifolds for $\varphi$ needs to be made by means of the local stable and unstable manifold theorem (see \cite{moham99}, Theorem 3.1). Therein it is considered a homeomorphism $H(\omega )$ which establishes the equivalence of the geometrical structures arising for both sets of manifolds, and as a consequence the manifolds for $\varphi$ inherit the same dynamics as the ones for $\Phi (t, \omega )$ but only in a neighborhood of the stationary orbit. In this sense the existence of such manifolds for a nonlinear RDS $\varphi$ is only ensured locally. Anyway this result provides a very good approximation to the stochastic dynamics of a system, and enables us to discuss the different patterns of behavior of the solutions in the following examples.}. In this direction there is a recent result (\cite{cheng16}, Theorem 2.1) which ensures the equivalence in the dynamics of both kinds of equations when the noisy term $\sigma$ is additive (i.e., $\sigma$ does not depend on the solution $X_{t}$). Although this was done by means of the most probable phase portrait, a technique that closely resembles the ordinary phase space for deterministic systems, this fact might indicate that such analogy in the dynamics cannot be achieved when the noise does depend explicitly on the solution $X_{t}$. Actually any additive noise affects all the particles together at the same magnitude.

Anyway the noisy saddle serves to establish an analogy to the dynamics with the deterministic saddle. One of its features is the contrast between the growth of the components $X_{t}$ and $Y_{t}$, which mainly have a positive and negative exponential growth respectively. We will see that this is graphically captured when applying the stochastic Lagrangian descriptors method to the SDE (\ref{noisy_saddle}) over a domain of the stationary orbit. Moreover when representing the stochastic Lagrangian descriptor values for the noisy saddle, one can observe that the lowest values are precisely located on the manifolds $\mathcal{S}(\omega )$ and $\mathcal{U}(\omega )$. These are manifested as  sharp features indicating a rapid change of the values that the stochastic Lagrangian descriptor assumes. This geometrical structure formed by ``local minimums" has a very marked crossed form and it is straightforward to think that the stationary orbit is located at the intersection of the two cross-sections. These statements are supported afterwards by numerical simulations and analytical results.

However there persists an open question about how reliable the stochastic Lagrangian descriptors method is when trying to depict the phase space of an arbitrary stochastic differential equation. This is the main issue concerning this method and has only been partially reported for deterministic dynamical systems in previous articles (\cite{MCWGM15} and \cite{LBGWM15}). In this last paper it is analytically proven the efficacy of this method for autonomous and non-autonomous Hamiltonian systems. The theoretical idea that supports this assertion is the discontinuity of the transversal derivative of the Lagrangian descriptor function over the manifolds of the corresponding hyperbolic trajectory. Following this idea, these ``singular features" arising on the manifolds of a hyperbolic trajectory for a deterministic Hamiltonian system motivates  us to study whether the ``abrupt changes" on the stochastic Lagrangian descriptor function represent the location of the manifolds of a stationary orbit or not. Another related question is to determine the size of the random term $\sigma dW$ in relation to  its influence on  the phase space of the deterministic equation $dX = b dt$. The next sections in this paper will be dedicated to addressing  these  issues by considering concrete examples of SDEs.

\section{The stochastic Lagrangian descriptor}
\label{sec:SLD}

Lagrangian descriptors have been shown to be a useful technique for revealing phase space structures in dynamical systems. The original Lagrangian 
descriptor, introduced in \cite{cnsns} for continuous dynamical systems, corresponded to the Euclidean arc length of a 
trajectory  over a time interval (backwards and forwards). In particular, we  consider the general  time-dependent vector field

\begin{equation}
\displaystyle{\frac{dx}{dt} = v(x,t), \quad x \in \mathbb{R}^n, t \in \mathbb{R},}
\end{equation}

\noindent
where $v(x,t) \in C^r (r \geq 1)$ in $x$ and continuous in time. For any solution $x(t) \equiv x(t, t_0, x_0)$  with initial condition  $x_0 \equiv x(t_0, t_0, x_0) \in 
\mathbb{R}^n$ and a fixed integration time $\tau$, the Lagrangian descriptor was initially defined as

\begin{equation}
\displaystyle{M(x_0, t_0, \tau) = \int^{t_0 + \tau}_{t_0 - \tau} \norm{\dot{x}(x_0, t_0, t)}dt},
\end{equation}

\noindent
where $\norm{\cdot }$ is the Euclidean norm. Afterwards in \cite{carlos}, the Lagrangian descriptor was adapted to the 
discrete dynamical systems setting but including a small change in the chosen norm ($p$-norm). Let $\{x_i\}^{i = N}_{i = 
-N}$, $N \in \mathbb{N}$ denote an orbit of $(2N + 1)$ nodes long and $x_i \in \mathbb{R}^n$. Considering the space of 
orbits as a sequence space, the discrete Lagrangian descriptor was defined in terms of the $\ell^p$-norm of an orbit as follows:

\begin{equation}
\displaystyle{MD_p(x_0, N) = \sum^{N-1}_{i = -N}\norm{x_{i+1} - x_i}_p, \quad p \leq 1.}
\end{equation}

\noindent
This alternative definition allows us to prove formally the non-differentiability of the $MD_p$ function through points that 
belong to invariant manifolds of a hyperbolic trajectory. This fact implies a better visualization of the invariant 
manifolds as they are detected over areas where the $MD_p$ function presents abrupt changes in its values.

Later on this last definition was adapted to the continuous time case in \cite{LBGWM15}. For any initial condition $x_0 = 
x(t_0) \in \mathbb{R}^n$ and any given time interval $[t_0 - \tau, t_0 + \tau]$, the $M$ function is redefined as

\begin{equation}
\displaystyle{M_p(x_0, t_0, \tau) = \int^{t_0 + \tau}_{t_0 - \tau} \norm{\dot{x}(x_0,t)}_pdt}
\end{equation}

\noindent
where $p$ is chosen from the interval $(0, 1]$.

Now we extend  these ideas to the context of stochastic differential equations. 
For this purpose we consider  a general SDE of the form:

\begin{equation}
dX_t = b(X_t, t)dt + \sigma(X_t, t)dW_t, \quad X_{t_0} = x_0,
\label{eq:stochastic_lagrangian_system}
\end{equation}

\noindent
where $X_t$ denotes the solution of the system, $b(\cdot)$ and $\sigma(\cdot)$ are Lipschitz functions which ensure 
uniqueness of solutions and $W_t$ is a two-sided Wiener process. Henceforth, we make use of the following notation

\begin{equation}
X_{t_j} := X_{t_0+j\Delta t},
\end{equation}

\noindent
for a given $\Delta t>0$ small enough and $j=-N,\cdots ,-1,0,1, \cdots ,N$.

\begin{definition}
The stochastic Lagrangian descriptor evaluated for SDE \eqref{eq:stochastic_lagrangian_system} with general solution 
$\textbf{X}_{t}(\omega )$ is given by
\begin{equation}
\displaystyle{MS_p(\textbf{x}_0, t_0, \tau, \omega) = \sum^{N-1}_{i = -N}\norm{\textbf{X}_{t_{i+1}} - 
\textbf{X}_{t_i}}_p,}
\label{eq:MS}
\end{equation}
where $\lbrace \textbf{X}_{t_{j}} \rbrace_{j=-N}^{N}$ is a discretization of the solution such that 
$\textbf{X}_{t_{-N}} = \textbf{X}_{-\tau}$, $\textbf{X}_{t_N} = \textbf{X}_{\tau}$, $\textbf{X}_{t_0} = \textbf{x}_{0}$, 
for a given $\omega \in \Omega$.
\end{definition}

\noindent
Obviously every output of the $MS_p$ function highly depends on the experiment $\omega \in \Omega$ where $\Omega$ is 
the probability space that includes all the possible outcomes of a given phenomena. Therefore in order to analyze the 
homogeneity of a given set of outputs, we consider a sequence of results of the $MS_p$ function for the same 
stochastic equation \eqref{eq:stochastic_lagrangian_system}: $MS_p(\cdot, \omega_1)$, $MS_p(\cdot, \omega_2)$, 
$\cdots$, $MS_p(\cdot, \omega_M)$. It is feasible that the following relation holds

\begin{equation}
d(MS_p(\cdot, \omega_i), MS_p(\cdot, \omega_j)) < \delta, \quad \text{for all } i,j,
\label{eq:deterministic_tol}
\end{equation}

\noindent
where $d$ is a metric that measures the similarity between two matrices (for instance $\norm{A-B}_F = 
\sqrt{Tr((A-B)\cdot (A-B)^T)}$ - Frobenius norm) and $\delta$ a positive tolerance. Nevertheless for general stochastic 
differential equations, expression (\eqref{eq:deterministic_tol}) does not usually hold. Alternatively if the elements 
of the sequence of matrices $MS_p(\cdot, \omega_1)$, $MS_p(\cdot, \omega_2)$, $\cdots$, $MS_p(\cdot, \omega_M)$ do not 
have much similarity to each other, it may be of more use to define the mean of the outputs

\begin{equation}
\displaystyle{\mathbb{E} \left[ MS_p(\cdot, \omega) \right] = \left(
\frac{MS_p(\cdot, \omega_1) + MS_p(\cdot, \omega_2) + \cdots +
MS_p(\cdot, \omega_M)}{M}\right) ,}
\label{eq:mean_MSp_value}
\end{equation}

\noindent
for a sufficiently large number of experiments $M$. Since the solution of a SDE is affected by the noise term, the phase 
portrait of the studied SDE for an arbitrary $\omega$ becomes unpredictable and one can only refer to places where 
hyperbolic trajectories and invariant manifolds are likely located. This way of understanding the geometry of transport 
for SDEs is similar in spirit as the one explained in \cite{banisch16} where the authors provide an alternative method 
for revealing coherent sets from Lagrangian trajectory data.

\section{Numerical simulation of the stochastic Lagrangian descriptor}
\label{sec:num}

In this section we describe the method  numerically solving the SDEs that we use  throughout this article.
Consider a general $n$-dimensional SDE of the form

\begin{equation}
dX^j_t = b^j(X_t, t)dt + \sum^m_{k=1}\sigma^j_k(X_t, t)dW^k_t, \quad X_{t_0} = x_0 \in \mathbb{R}^n, \quad j=1,\cdots ,n
\end{equation}

\noindent
where $X_t = (X^1_t, \cdots , X^n_t)$ and $W^1_t, \cdots, W^m_t$ are $m$ independent Wiener processes.
If the time step $\Delta t$ is firstly fixed, then the temporal grid $t_p = t_0 + p\Delta t$ ($p \in \mathbb{Z}$) is 
already determined and we arrive to the difference equation

\begin{equation}
X^j_{t+\Delta t} = X^j_t + b^j(X_t, t)\Delta t + \sum^m_{k=1} \sigma^j_k(X_t, t)dW^k_t.
\end{equation}

\noindent
This scheme is referred to as  the Euler-Marayuma method for solving a single path of the SDE. If the stochastic part is removed 
from the equation, then we obtain the classical Euler method. Suppose $X_{t_p}$ is the solution of the SDE and 
$\tilde{X}_{t_p}$ its numerical approximation at any time $t_p$. Since both of them are random variables, the accuracy 
of the method must be determined in probabilistic terms. With this aim, the following definition is introduced.

\begin{definition}
A stochastic numerical method has an  order of convergence equal to $\gamma$ if there exists a constant $C>0$ such that

\begin{equation}
\mathbb{E} \left[ X_{t_p} - \tilde{X}_{t_p} \right] \leq C \Delta t^{\gamma},
\end{equation}
for any arbitrary $t_p = t_0 + p\Delta t$ and $\Delta t$ small enough.
\end{definition}

\noindent
Indeed, the Euler-Maruyama method has an order of convergence equal to $1/2$ (see \cite{kloeden92} for further details).

\section{The noisy saddle}

The noisy saddle is a fundamental benchmark for assessing numerical methods for revealing phase space structures.
Its main value is in the simplicity of the expressions taken by the components of the stationary orbit and its  corresponding stable and unstable manifolds. From these one clearly observes the exponential separation rates between particles passing 
near the manifolds. Now for the stochastic differential equations

\begin{equation}
\label{eq:general_noisy}
\begin{cases}
  dX_t = a_1X_t dt + b_1dW^1_t \\
  dY_t = -a_2Y_t dt + b_2dW^2_t
\end{cases}
\end{equation}

\noindent
it is straightforward to check that the only stationary orbit takes the expression

\begin{equation}
  \widetilde{X}(\omega ) = \left( \tilde{x}(\omega ), \tilde{y}(\omega ) \right) = \left( 
\displaystyle{-\int_{0}^{\infty}e^{-a_{1}s} b_1dW^1_t(\omega )} , \displaystyle{\int_{-\infty}^{0}e^{b_{1}s} 
b_2dW^2_t(\omega )} \right)
\end{equation}

\noindent
where $a_{1},a_{2},b_{1},b_{2} \in \mathbb{R}$ are constants and $a_{1},a_{2}>0$. Its corresponding stable and unstable 
manifolds are

\begin{equation}
\mathcal{S}(\omega ) = \lbrace (x,y) \in \mathbb{R}^{2} : x = \tilde{x}(\omega ) \rbrace , \quad
\mathcal{U}(\omega ) = \lbrace (x,y) \in \mathbb{R}^{2} : y = \tilde{y}(\omega ) \rbrace .
\end{equation}

\noindent
These play a very relevant role as dynamical barriers for the particles tracked by the RDS, which is generated by the 
SDE (\ref{eq:general_noisy}). This fact has been justified in the previous section, but can be analytically demonstrated  when 
computing the stochastic Lagrangian descriptor (\ref{eq:MS}) for the solution of the noisy saddle.

According to the notation used for the definition of SLD

$$\displaystyle{MS_p(\textbf{x}_0, t_0, \tau, \omega) = \sum^{N-1}_{i = -N}\norm{\textbf{X}_{t_{i+1}} - 
\textbf{X}_{t_i}}_p,}$$

\noindent
at which the components of the solution satisfy the initial conditions $\textbf{X}_{t_{0}}= \left( X_{t_{0}},Y_{t_{0}} 
\right) = (x_{0},y_{0}) = \textbf{x}_{0}$, these take the expressions

\begin{equation}
\label{general_noisy_saddle_solutions}
X_{t} = e^{a_{1}t} \left( x_{0} + \int_{0}^{t}e^{-a_{1}s}b_{1}dW_{s}^{1} \right) \quad , \quad Y_{t} = e^{-a_{2}t} 
\left( y_{0} + \int_{0}^{t}e^{a_{2}s}b_{2}dW_{s}^{2}(\omega ) \right)
\end{equation}

\noindent
and the temporal nodes satisfy the rule $t_{i} = t_{0} + i\Delta t$ with $t_{0}$ and $\Delta t$ already given. Now it is 
possible to compute analytically the increments $\norm{\textbf{X}_{t_{i+1}} - \textbf{X}_{t_i}}_p = \left| X_{t_{i+1}} - 
X_{t_i} \right|^{p} + \left| Y_{t_{i+1}} - Y_{t_i} \right|^{p}$:

$$\left| X_{t_{i+1}} - X_{t_i} \right|^{p} = \left| e^{a_{1}t_{i+1}} \left( x_{0} + 
\int_{0}^{t_{i+1}}e^{-a_{1}s}b_{1}dW_{s}^{1} \right) - e^{a_{1}t_{i}} \left( x_{0} + 
\int_{0}^{t_{i}}e^{-a_{1}s}b_{1}dW_{s}^{1} \right) \right|^{p}$$

$$= \left| e^{a_{1}t_{i}}\left( e^{a_{1}\Delta t} - 1 \right) \left[ x_{0} + \int_{0}^{t_{i}}e^{-a_{1}s}b_{1}dW_{s}^{1} 
\right] + e^{a_{1}(t_{i}+\Delta t)} \int_{t_{i}}^{t_{i}+\Delta t}e^{-a_{1}s}b_{1}dW_{s}^{1} \right|^{p}$$
$$= \left| e^{a_{1}t_{i}}\left( e^{a_{1}\Delta t} - 1 \right) \left[ x_{0} + 
\int_{0}^{t_{i}}e^{-a_{1}s}b_{1}dW_{s}^{1} \right] + e^{a_{1}\Delta t}b_{1}dW_{t_{i}}^{1} \right|^{p} \text{(by applying 
It\^{o} formula (\ref{Ito})).}$$

\noindent
Moreover for large values of $t_{i}$ such that $e^{a_{1}t_{i}} \gg e^{a_{1}\Delta t}$ and taking into account that 
$dW_{t_{i}}^{1}$ is finite almost surely, we can consider the following approximation

\begin{equation}
\label{x_increments}
\left| X_{t_{i+1}} - X_{t_i} \right|^{p} \hspace{0.1cm} \approx \hspace{0.1cm} e^{a_{1}t_{i}\cdot p}\left| 
e^{a_{1}\Delta t} - 1 \right|^{p} \left| x_{0} + \int_{0}^{t_{i}}e^{-a_{1}s}b_{1}dW_{s}^{1} \right|^{p}.
\end{equation}

\noindent
By following these arguments, one can get an analogous result for the second component $Y_{t}$:

$$\left| Y_{t_{i+1}} - Y_{t_i} \right|^{p} = \left| e^{-a_{2}t_{i}}\left( e^{-a_{2}\Delta t} - 1 \right) \left[ 
y_{0} + \int_{0}^{t_{i}}e^{a_{2}s}b_{2}dW_{s}^{2} \right] + e^{-a_{2}\Delta t}b_{2}dW_{t_{i}}^{2} \right|^{p},$$

\noindent
which for small values of $t_{i}$, such that $e^{-a_{2}t_{i}} \gg e^{-a_{2}\Delta t}$, this approximation can be further 
simplified as follows

\begin{equation}
\label{y_increments}
\left| Y_{t_{i+1}} - Y_{t_i} \right|^{p} \hspace{0.1cm} \approx \hspace{0.1cm} e^{-a_{2}t_{i}\cdot p}\left| 
e^{-a_{2}\Delta t} - 1 \right|^{p} \left| y_{0} + \int_{0}^{t_{i}}e^{a_{2}s}b_{2}dW_{s}^{2} \right|^{p}.
\end{equation}

Once the analytic expression of the SLD applied to the noisy saddle (\ref{eq:general_noisy}) is known, it can be proved 
that the stable and unstable manifolds of the stationary orbit are manifested as singularities of the SLD function over any given domain of 
initial conditions containing the stationary orbit. This fact implies that the SLD method realizes a procedure to detect 
these geometrical objects and, consequently, provides a phase portrait representation of the dynamics generated by the 
noisy saddle. In the same way as described in \cite{LBGWM15}, we refer to singularities as points of the domain of spatial  initial conditions where the derivative of the SLD is not defined. The paradigm example of the mathematical manifestation of singularities of the LD on stable and unstable manifolds of hyperbolic trajectories is provided by the scalar function $|\cdot |^{p}$ with 
$p \in (0,1]$.  This function is  singular, alternatively non-differentiable, at those points where its argument is  zero. Graphically this feature is observed as sharp changes in the representation of the SLD values, where the 
contour lines  concentrate in a very narrow space. 
\newline
\\
In this particular example we are able to explicitly identify within the expression of the SLD  the terms that are largest in 
magnitude. In other words, we are able to identify the terms whose particular singularities determine the non-differentiability of the 
entire sum\footnote{Note that the differentiability of the SLD is analyzed with respect to the 
components of the initial conditions $(x_{0},y_{0})$, as the experiment $\omega \in \Omega$ and the starting point 
$t_{0}$ are previously fixed.}. This is better understandable if the expression of the SLD is divided into two sums

$$\displaystyle{MS_p(\textbf{x}_0, t_0, \tau, \omega) = \sum^{N-1}_{i = -N}\norm{\textbf{X}_{t_{i+1}} - 
\textbf{X}_{t_i}}_p = \sum^{N-1}_{i = -N}\left| X_{t_{i+1}} - X_{t_i} \right|^{p}} + \sum^{N-1}_{i = -N}\left| 
Y_{t_{i+1}} - Y_{t_i} \right|^{p} .$$

\noindent
The highest order term within the first sum is $\left| X_{t_{N}} - X_{t_{N-1}} \right|^{p} = \left| X_{\tau } - X_{\tau 
- \Delta t} \right|^{p}$, which according to (\ref{x_increments})  is approximated by

\begin{equation}
\label{higher_order_x}
\left| X_{\tau } - X_{\tau - \Delta t} \right|^{p} \hspace{0.1cm} \approx \hspace{0.1cm} e^{a_{1}(\tau - \Delta 
t)\cdot p}\left| e^{a_{1}\Delta t} - 1 \right|^{p} \left| x_{0} + \int_{0}^{\tau - \Delta t}e^{-a_{1}s}b_{1}dW_{s}^{1} 
\right|^{p} \quad \text{for enough large values of } \tau .
\end{equation}

\noindent
Similarly the highest order term within the second sum is $\left| Y_{t_{-N+1}} - Y_{t_{-N}} \right|^{p} = \left| 
Y_{-\tau +\Delta t} - Y_{-\tau } \right|^{p}$, approximated by

\begin{equation}
\label{higher_order_y}
\left| Y_{-\tau +\Delta t} - Y_{-\tau } \right|^{p} \hspace{0.1cm} 
\approx \hspace{0.1cm} e^{a_{2}\tau \cdot p}\left| e^{-a_{2}\Delta t} - 1 \right|^{p} \left| y_{0} - 
\int_{-\tau}^{0}e^{a_{2}s}b_{2}dW_{s}^{2} \right|^{p} \quad \text{for enough large values of } \tau .
\end{equation}

Consequently, it is evident that the sharper features will be located closed to the points where these two last 
quantities (\ref{higher_order_x}), (\ref{higher_order_y}) are zero. In other words where the initial condition 
$(x_{0},y_{0})$ satisfies one of the two following

$$x_{0} = - \int_{0}^{\tau - \Delta t}e^{-a_{1}s}b_{1}dW_{s}^{1} \quad \text{or} \quad y_{0} = 
\int_{-\tau}^{0}e^{a_{2}s}b_{2}dW_{s}^{2} \quad \text{for enough large values of } \tau .$$

\begin{figure}[htbp!]
\centering
{\includegraphics[scale = 0.5]{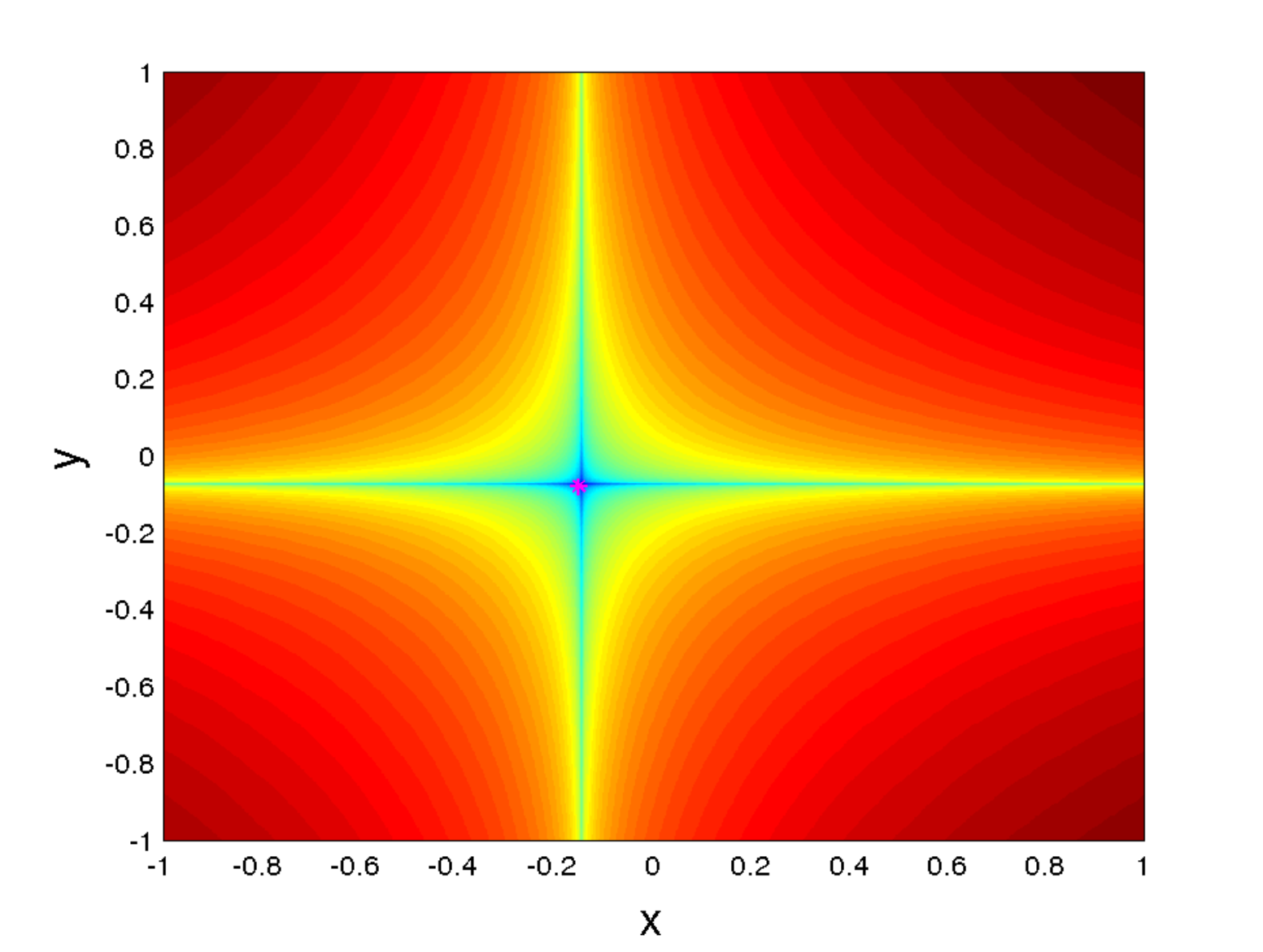}}
{\includegraphics[scale = 0.5]{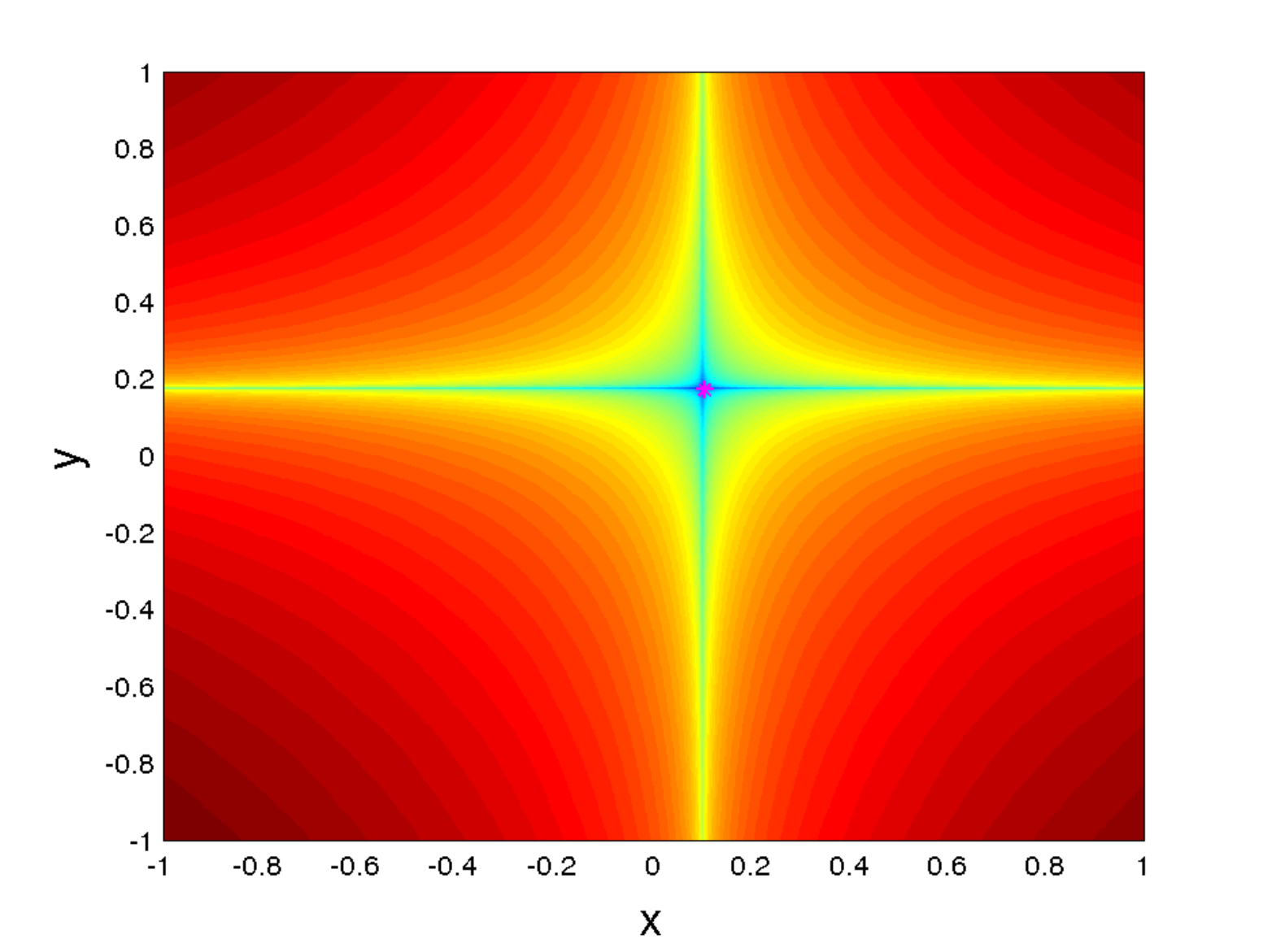}}
\caption{Two different experiments representing contours of $MS_p$ for $p=0.1$ and $\tau=15$. The contours of $MS_p$ are 
computed on a 1200$\times$ 1200 points grid of initial conditions and the time step for integration of the vector field 
is chosen to be $\Delta t= 0.05$. The magenta colored point corresponds to the location of the stationary orbit for each 
experiment. The chosen parameters are $a_1 = a_2 = b_2 = 1$ and $b_1 = -1$.}
\label{fig:saddle}
\end{figure}

\noindent
This statement is in agreement with the distinguished nature of the manifolds of the stationary orbit discussed in the
previous section. Note also that the two quantities for $x_{0}$ and $y_{0}$ converge to the coordinates of the 
stationary orbit $(\tilde{x}(\omega ),\tilde{y}(\omega ))$ with $\tau$ tending to infinity. These features are observed 
in Figure \ref{fig:saddle}, where the sharpness of the SLD representation highlights the location of the stable and unstable manifolds. 
The intersection of the two `'singular'' curves  represents the 
position of the stationary orbit $(\tilde{x}(\omega ),\tilde{y}(\omega ))$ for a given $\omega \in \Omega$. 
This fact is validated by the depiction of the stationary orbit, whose components have been computed separately from the 
SLD, and using the same output of the Wiener process.

\begin{remark}
Due to the properties of It\^{o} integrals, see for instance \cite{duan15}, the components of the stationary orbit 
satisfy

$$\mathbb{E} \left[ \tilde{x}(\omega ) \right] = \mathbb{E} \left[ - \int_{0}^{\infty}e^{-s}dW_{s}^{1} \right] = 0 
\quad , \quad \mathbb{E} \left[ \tilde{y}(\omega ) \right] = \mathbb{E} \left[ \int_{-\infty}^{0}e^{s}dW_{s}^{2} \right] 
= 0,$$

$$\mathbb{V} \left[ \tilde{x}(\omega ) \right] = \mathbb{E} \left[ \tilde{x}(\omega )^{2} \right] = \mathbb{E} \left[ 
\int_{0}^{\infty}e^{-2s}ds \right] = \frac{1}{2} \quad , \quad \mathbb{V} \left[ \tilde{y}(\omega ) \right] = \mathbb{E} 
\left[ \tilde{y}(\omega )^{2} \right] = \mathbb{E} \left[ \int_{-\infty}^{0}e^{2s}ds \right] = \frac{1}{2}.$$

\noindent
This means that the stationary orbit $ (\tilde{x}(\omega ),\tilde{y}(\omega ))$ is highly probable to be located closed to 
the origin of coordinates $(0,0)$, and  this feature is displayed in Figure \ref{fig:saddle}. This result gives 
more evidences and supports the similarities between the stochastic differential equation (\ref{noisy_saddle}) and the 
deterministic analogue system $\lbrace \dot{x}=x, \hspace{0.1cm} \dot{y}=-y \rbrace$ whose only fixed point is  
$(0,0)$.
\end{remark}

Therefore we can assert that the stochastic Lagrangian descriptor is a technique that provides a phase portrait 
representation of the dynamics generated by the noisy saddle equation (\ref{eq:general_noisy}).  Next we apply this same technique to further examples.

\section{Other examples}
\label{sec:examp}

\subsection{Stochastically forced Duffing equation}

Another classical problem is the Duffing oscillator. The deterministic version is given by

\begin{equation}
\label{eq:duffing_determ}
\ddot{x} = \alpha \dot{x} + \beta x + \gamma x^3 + \epsilon \cos(t).
\end{equation}

\noindent
If $\epsilon = 0$ the Duffing equation becomes time-independent, meanwhile for $\epsilon \neq 0$ the oscillator is a time-dependent system, where $\alpha$ is  the damping parameter, $\beta$ controls the 
rigidity of the system and $\gamma$ controls the size of the nonlinearity of the restoring force. The stochastically forced Duffing 
equation is studied in \cite{datta01} and can be written as follows:

\begin{equation}
\begin{cases}
  dX_t = \alpha Y_t, \\
  dY_t = (\beta X_t + \gamma X^3_t)dt + \epsilon dW_t.
\end{cases}
\end{equation}

\begin{figure*}[htbp!]
  \centering
  \subfigure[$MS_p$ contours for $\omega_1$]{\includegraphics[width=0.48\linewidth]{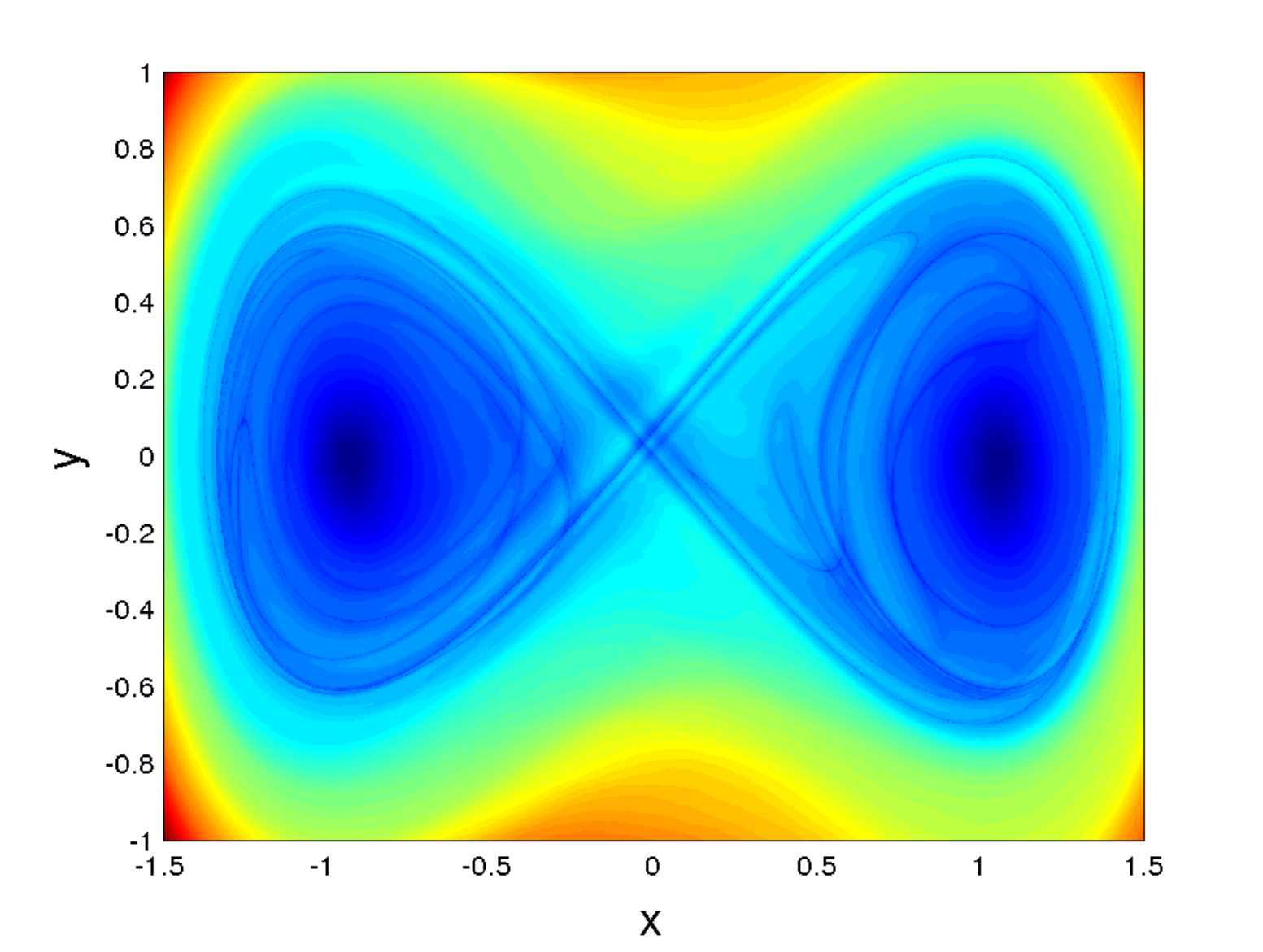}}
  \subfigure[$MS_p$ contours for $\omega_2$]{\includegraphics[width=0.48\linewidth]{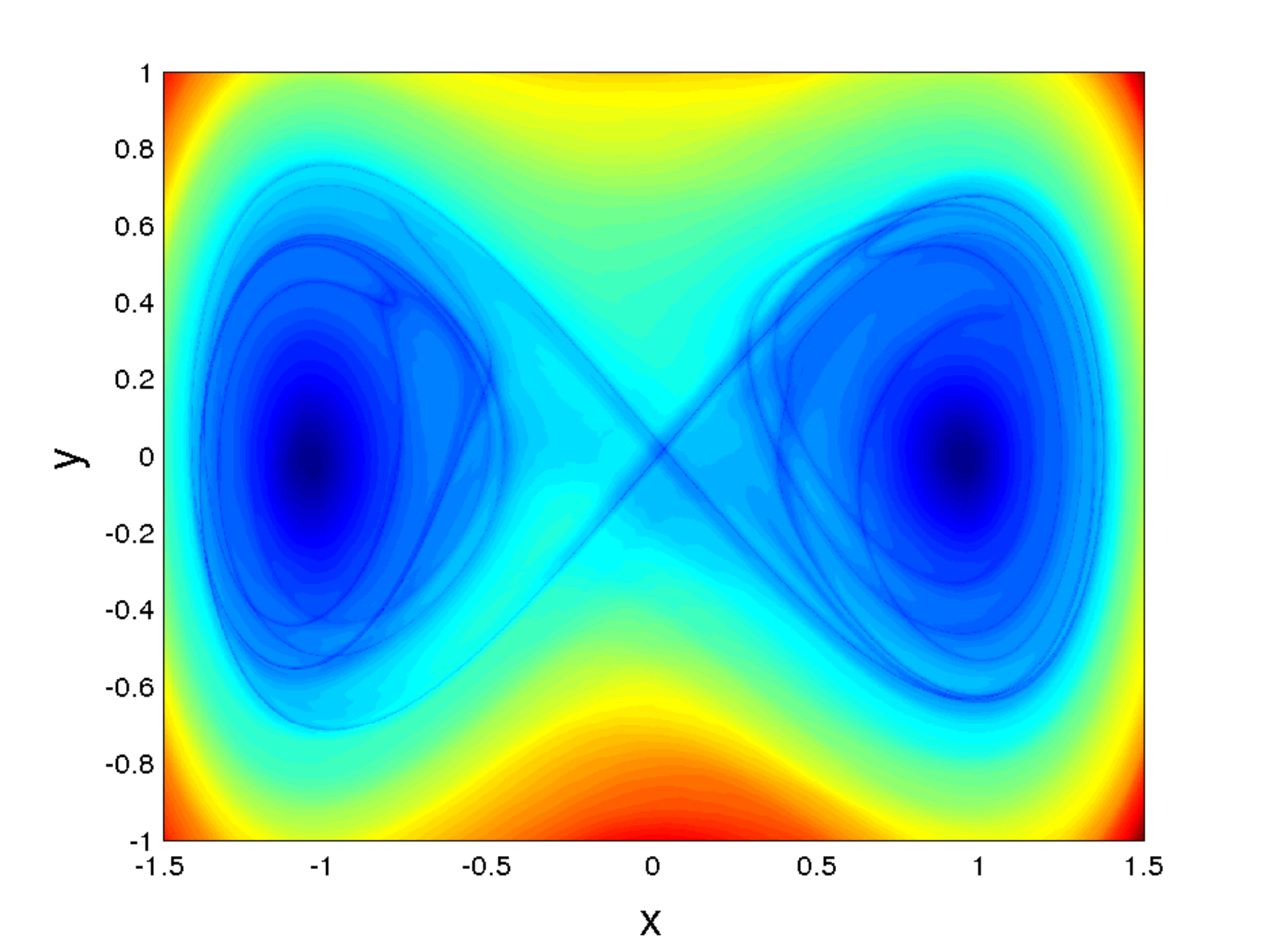}} \\
  \subfigure[$MS_p$ contours for $\omega_3$]{\includegraphics[width=0.48\linewidth]{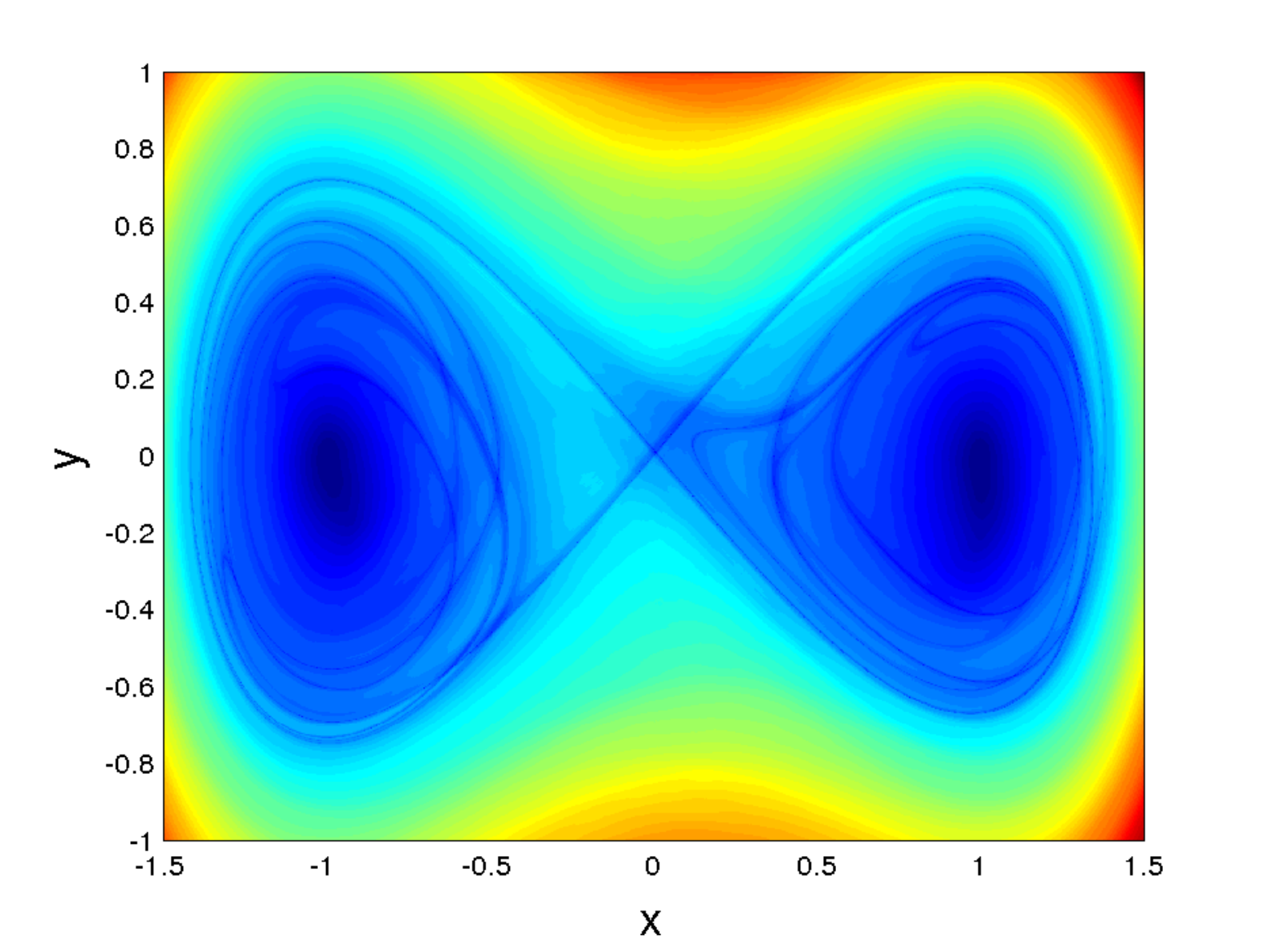}}
  \subfigure[$M_p$ for equation $\eqref{eq:duffing_determ}$]{\includegraphics[width=0.48\linewidth]{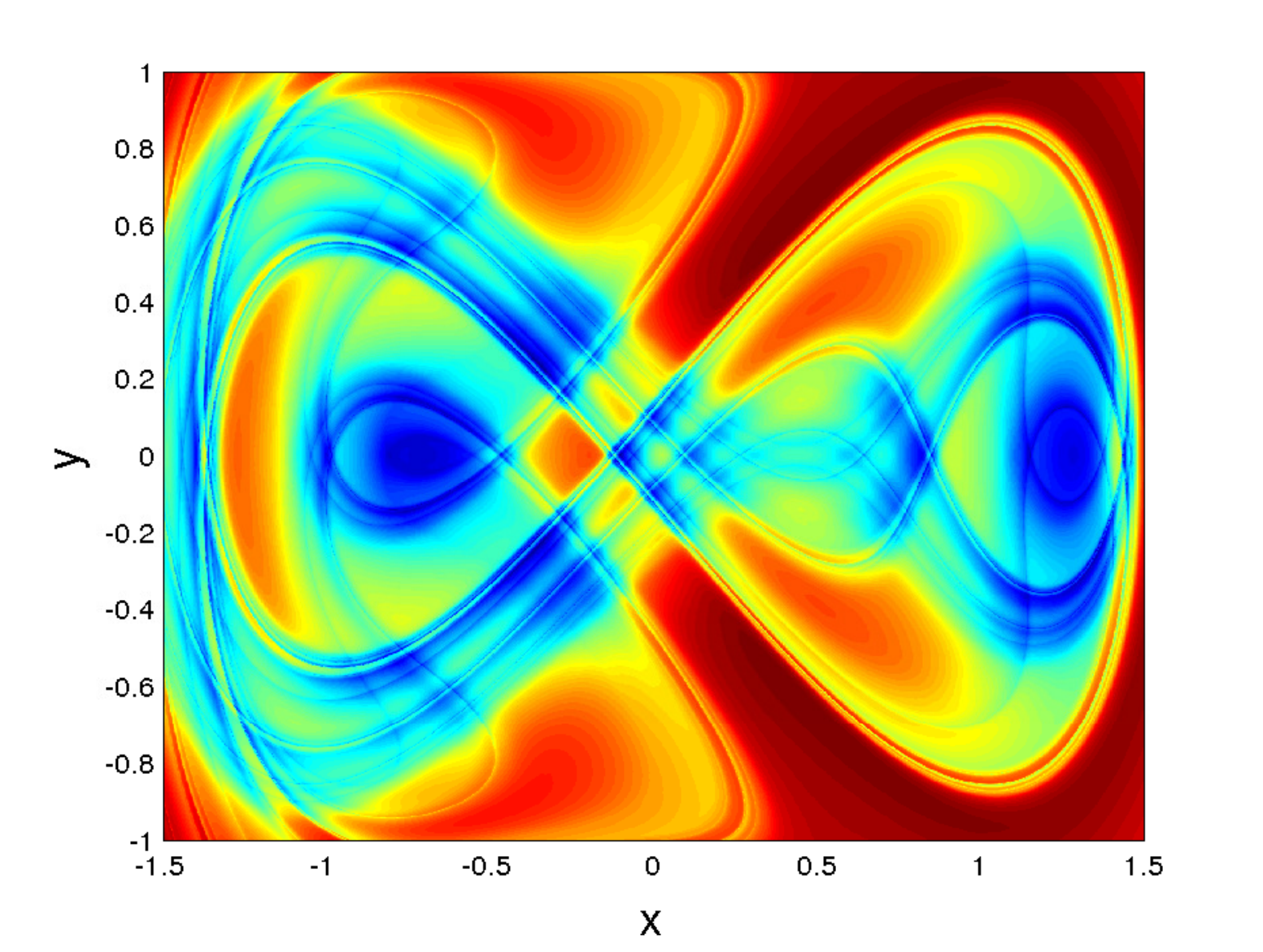}} \\
  \caption{a), b), c) Three different experiments representing $MS_p$ contours for $p=0.5$ over a grid of initial 
conditions. d) The last image corresponds to the $M_p$ function for equation \eqref{eq:duffing_determ} and $p=0.75$. All 
these pictures were computed for $\tau=15$, and over a $1200 \times 1200$ points grid. The time step for integration of 
the vector field was chosen to be $\Delta t = 0.05$.}
\label{fig:duffing}
\end{figure*}

\noindent
For our numerical  experiments we have selected $\alpha = \beta = 1$, $\gamma = -1$ and $\epsilon = 0.25$. The results of three 
different experiments (let say three different samples of $\omega_{1},\omega_{2},\omega_{3} \in \Omega$) are shown
in Figure \ref{fig:duffing}.

\subsection{Stochastically forced Double Gyre}

The double gyre is a standard benchmark velocity field in the study of the dynamical systems approach to Lagrangian transport. 
Recently it has been studied in the situation where the time-dependence is stochastic (\cite{hsieh12}). These equations are given by

\begin{equation}
\begin{cases}
  dX_t = \displaystyle{\left (-\pi A\sin\Big (\pi\frac{f(X_t,t)}{s}\Big ) \cos\Big (\pi\frac{Y_t}{s}\Big ) - \mu 
X_t\right )dt + \alpha dW^1_t }, \\ \\
  dY_t = \displaystyle{\left (\pi A\cos\Big (\pi\frac{f(X_t,t)}{s}\Big ) \sin\Big (\pi\frac{Y_t}{s}\Big ) 
\frac{\partial f}{\partial X_t} - \mu Y_t\right )dt + \alpha dW^2_t },
\end{cases}
\end{equation}

\noindent
where

\begin{equation}
f(X_t,t) = \epsilon \sin(\phi t + \psi)X_t^2 + (1 - 2\epsilon\sin(\phi t + \psi))X_t.
\end{equation}

\begin{figure}[htbp!]
\centering
{\includegraphics[scale = 0.5]{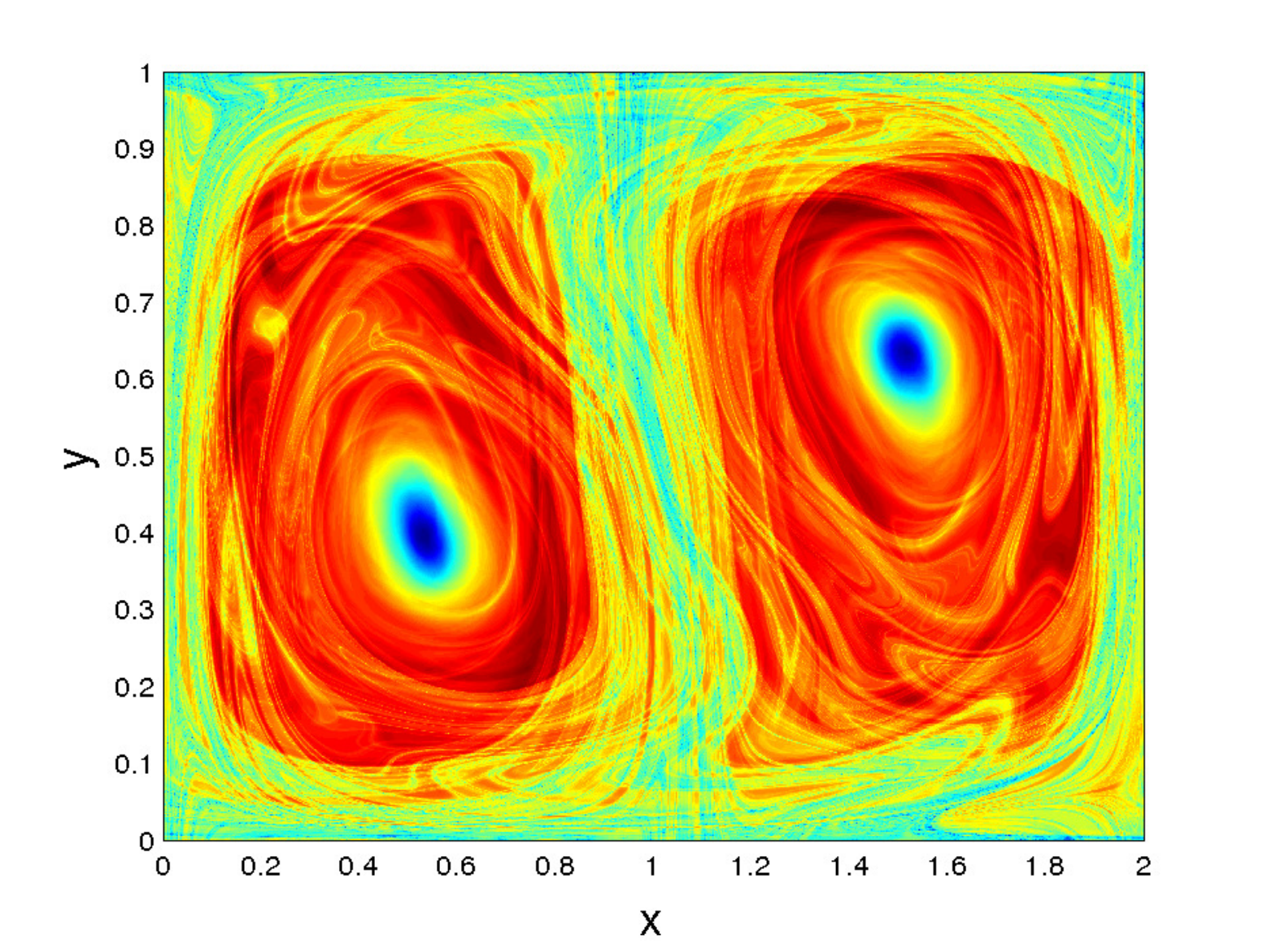}}
{\includegraphics[scale = 0.5]{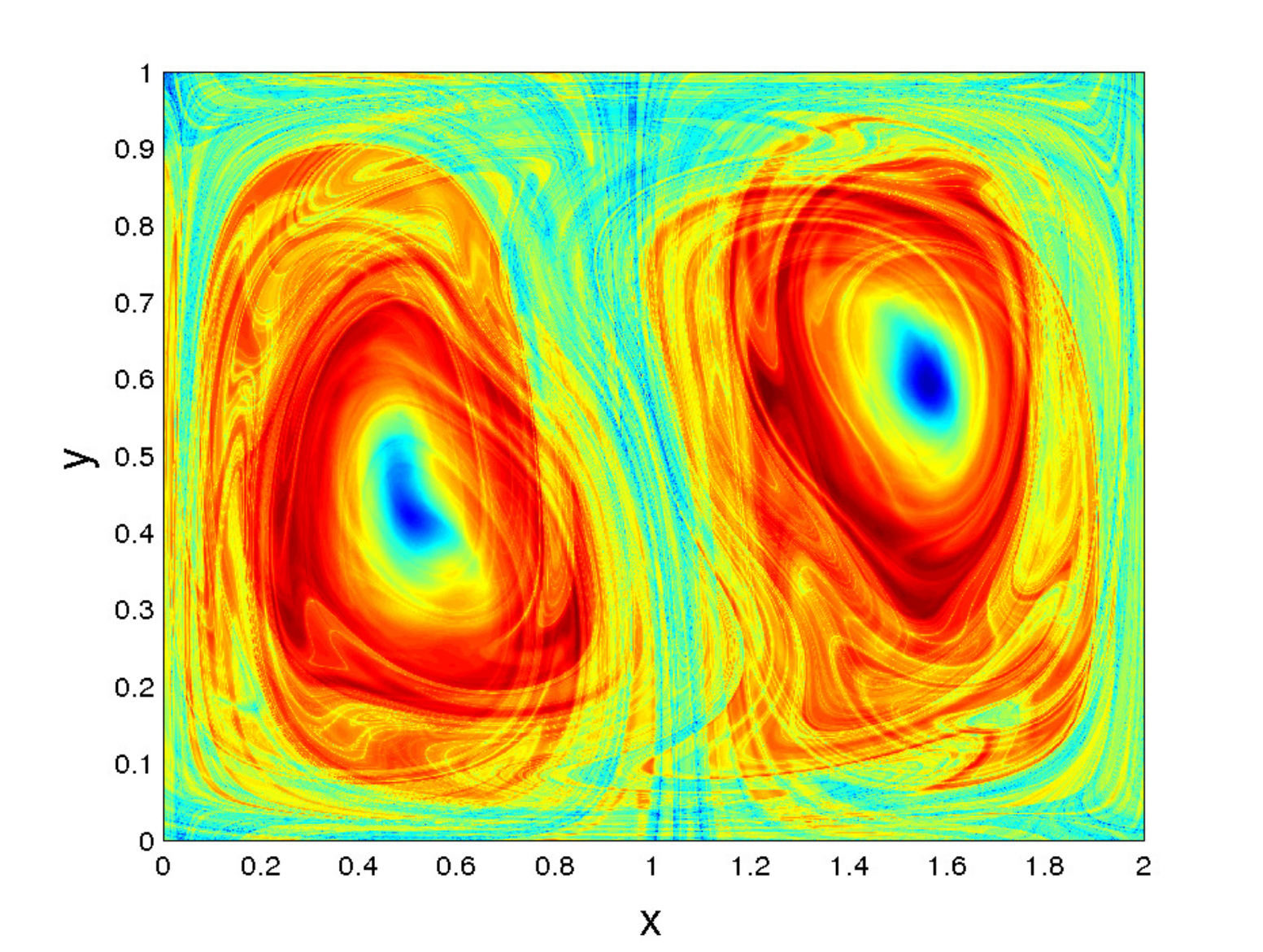}}
\caption{Two different experiments representing contours of $MS_p$ for $p=0.5$ and $\tau=15$. The contours of $MS_p$ are 
computed on a $1200\times 1200$ points grid of initial conditions and the time step for integration of the vector field 
is chosen to be $\Delta t = 0.05$.}
\label{fig:dgyre}
\end{figure}

\noindent
When $\epsilon = 0$ the double-gyre is time-independent, meanwhile if $\epsilon \neq 0$ the gyres force a periodic 
behavior in the $x$ direction. Among the parameters within the equation, $A$ models the amplitude of the velocity 
vectors, $\frac{\phi}{2\pi}$ gives the oscillation frequency, $\psi$ is the phase, $\mu$ determines the dissipation, $s$ 
scales the dimensions of the grid and $d W^i_t$ describes a Wiener process (stochastic white noise) with mean zero and standard 
deviation $\sigma = \Delta t$, while $\alpha$ is the amplitude of the noise. The results of two different experiments are 
observed in Figure \ref{fig:dgyre}. For such experiments we have used the following values for the parameters: $A = 
0.25$, $\phi = 2 \pi $, $\psi = 0$, $\mu = 0$, $s = 1$, $\alpha = 0.1$, $\epsilon = 0.25$.

\begin{figure}[htbp!]
\centering
\includegraphics[scale=0.8]{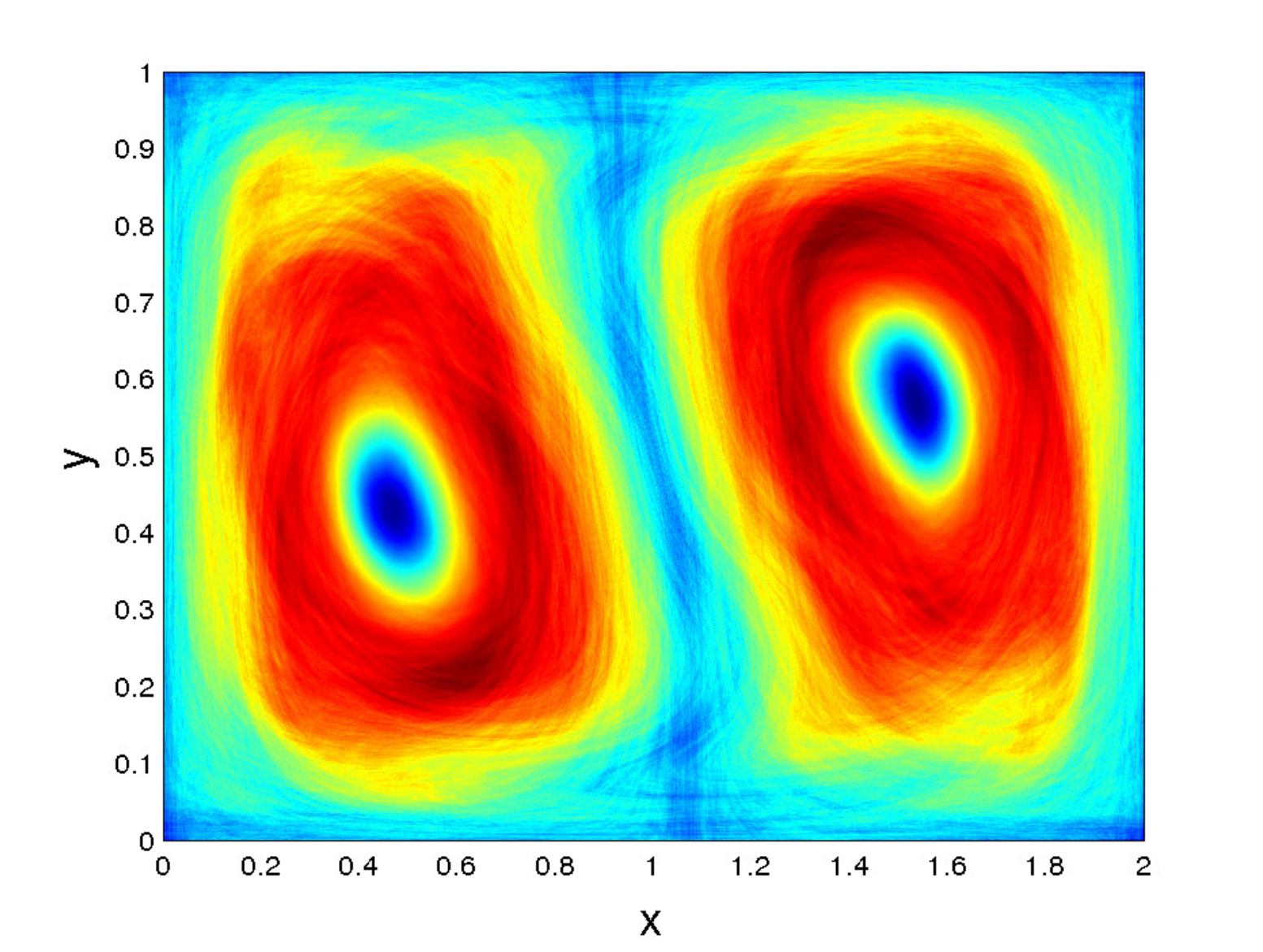}
\caption{Computation of $\mathbb{E} \left[ MS_p(\cdot, \omega) \right]$ for $p = 0.5$, $\tau = 15$ and 30 different 
experiments with a time step integration $\Delta t = 0.05$.}
\label{fig:30experiments}
\end{figure}

The outputs are quite different between these two experiments. Therefore we make  use of
\eqref{eq:mean_MSp_value}. In Figure \ref{fig:30experiments} we can observe the values of $\mathbb{E}\left[ MSp \right]$ when considering 30 different experiments. This is the expected phase space structure: two gyre centers are detected 
near the points $(0.5,0.5)$ and $(1.5,0.5)$ in dark blue, while around the line $\lbrace x=1 \rbrace$ we observe bundles 
of invariant manifolds in light blue. Thus for almost every $\omega \in \Omega$ two centers are depicted close to the 
middle of each gyre and by separating each gyre, the invariant manifolds can be interpreted as transport barriers.

\par
\begin{figure*}[htbp!]
  \centering
  \subfigure{\includegraphics[width=0.48\linewidth]{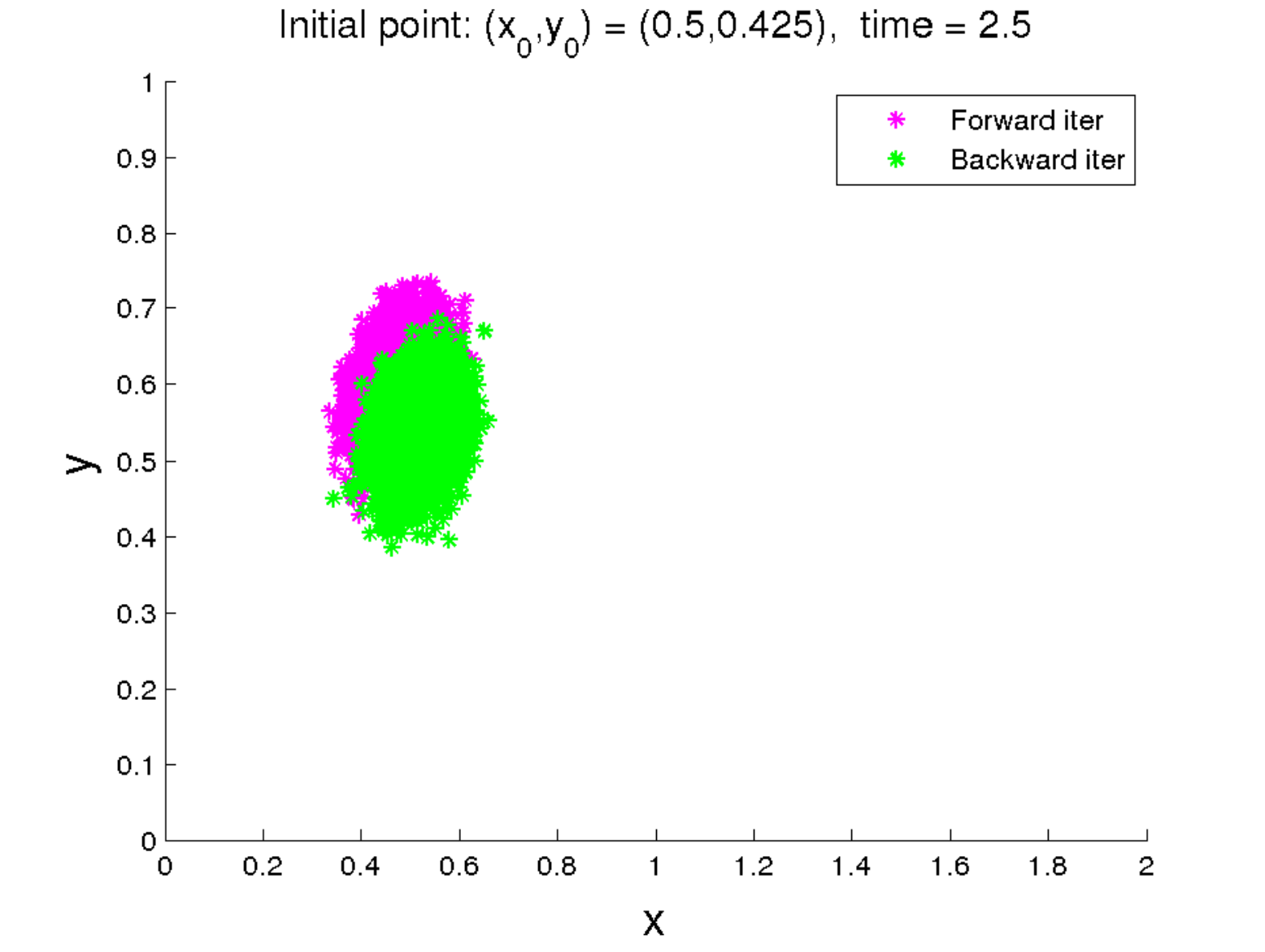}}
  \subfigure{\includegraphics[width=0.48\linewidth]{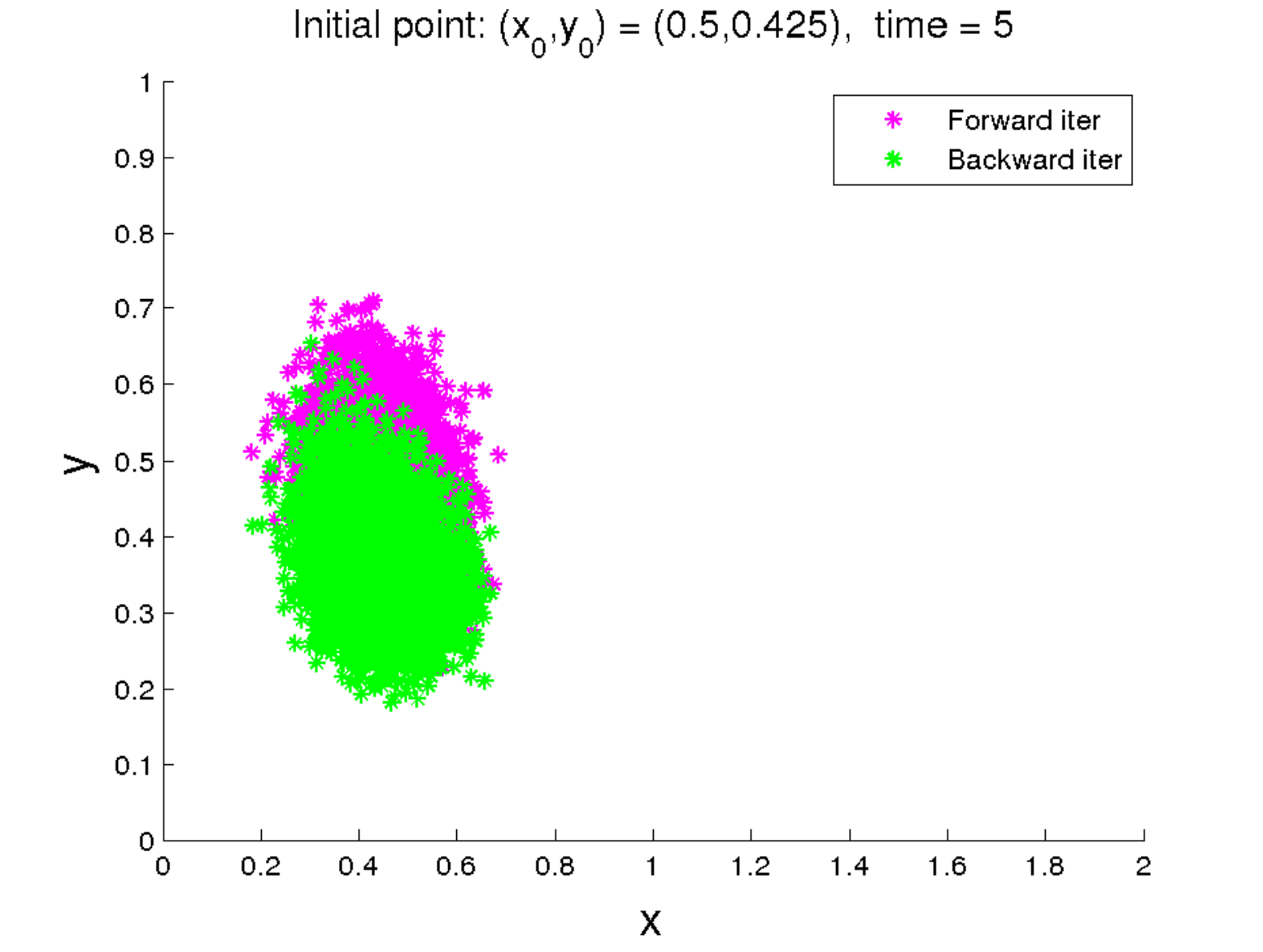}} \\
  \subfigure{\includegraphics[width=0.48\linewidth]{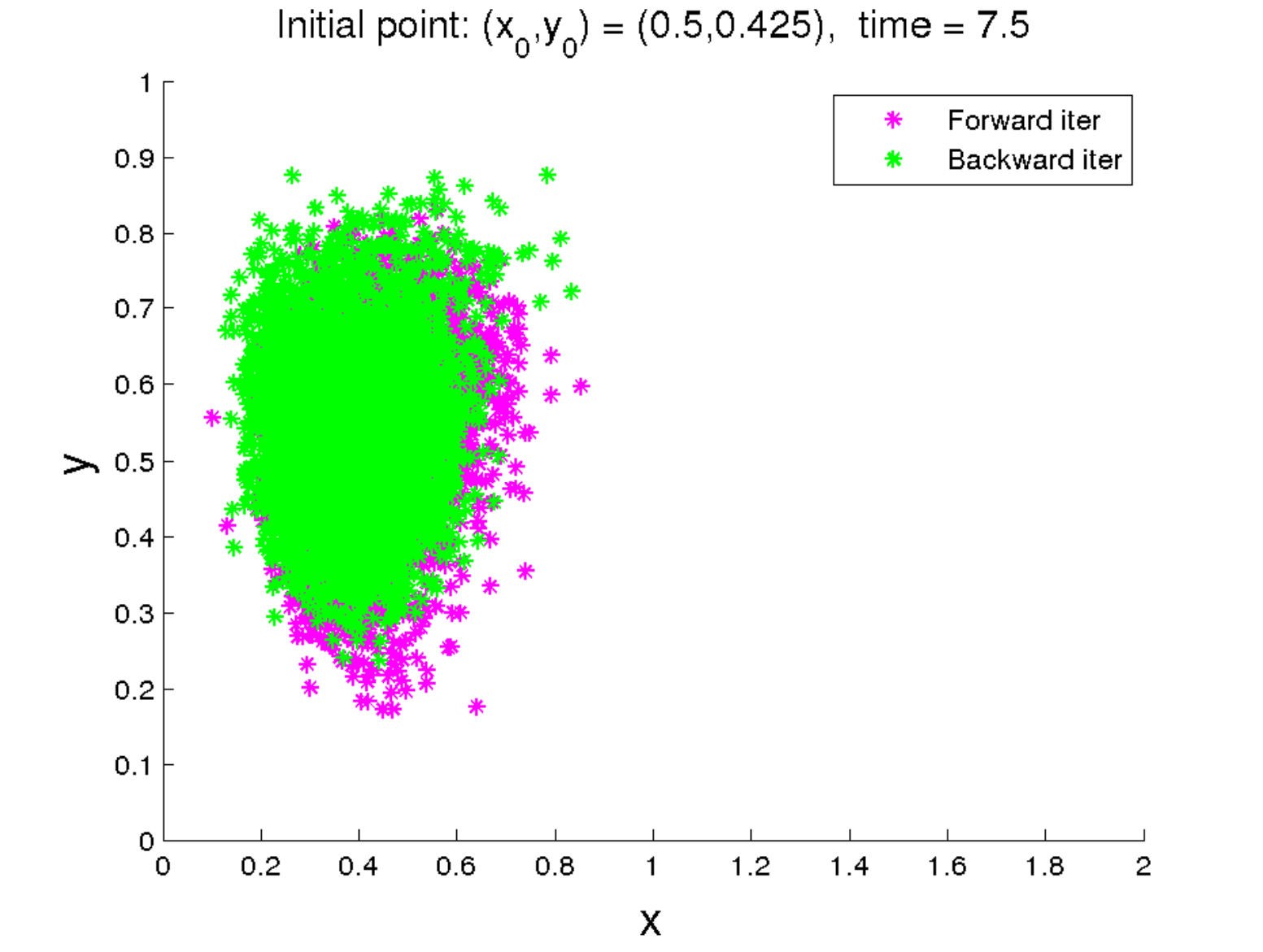}}
  \subfigure{\includegraphics[width=0.48\linewidth]{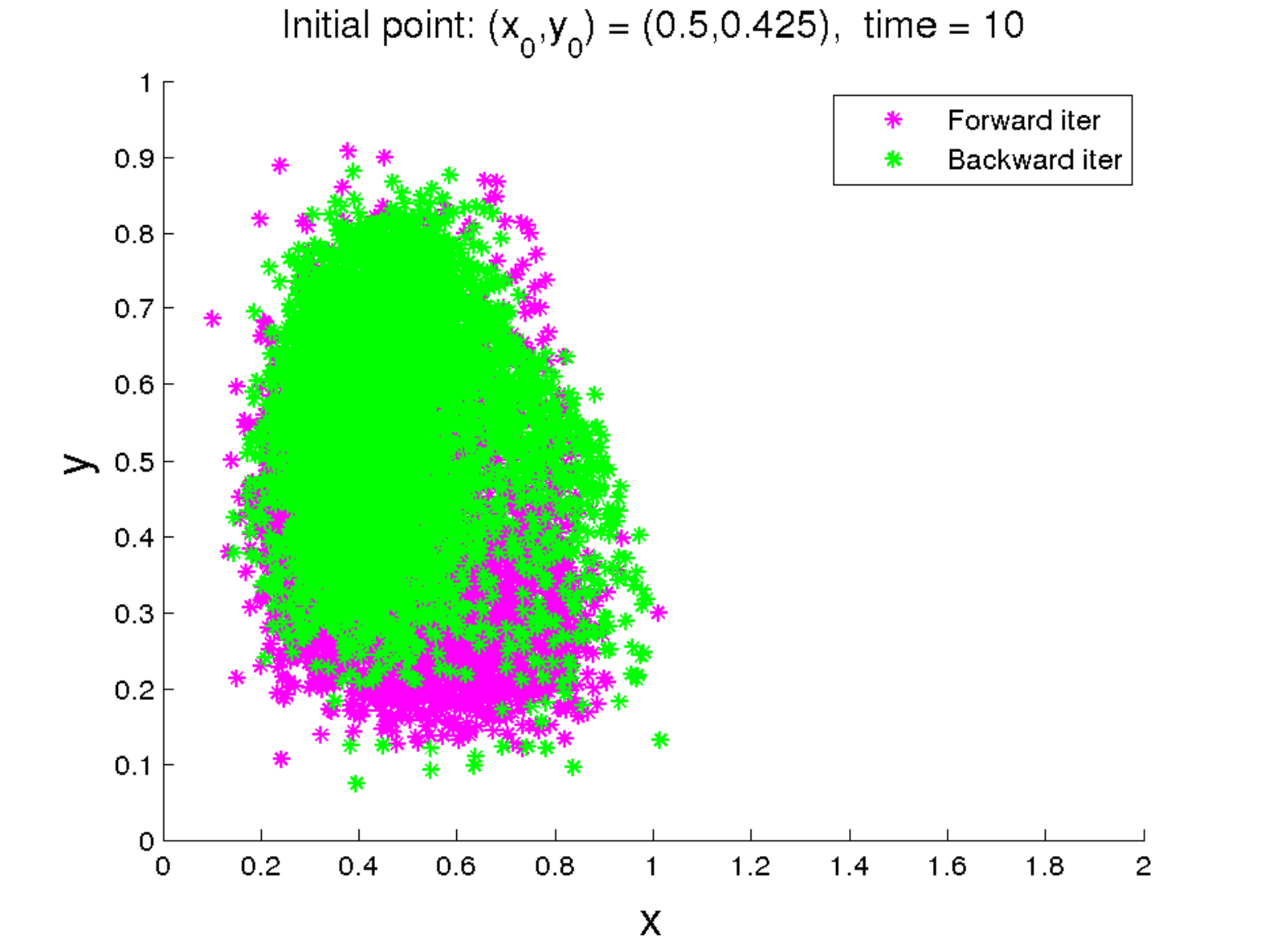}} \\
  \subfigure{\includegraphics[width=0.48\linewidth]{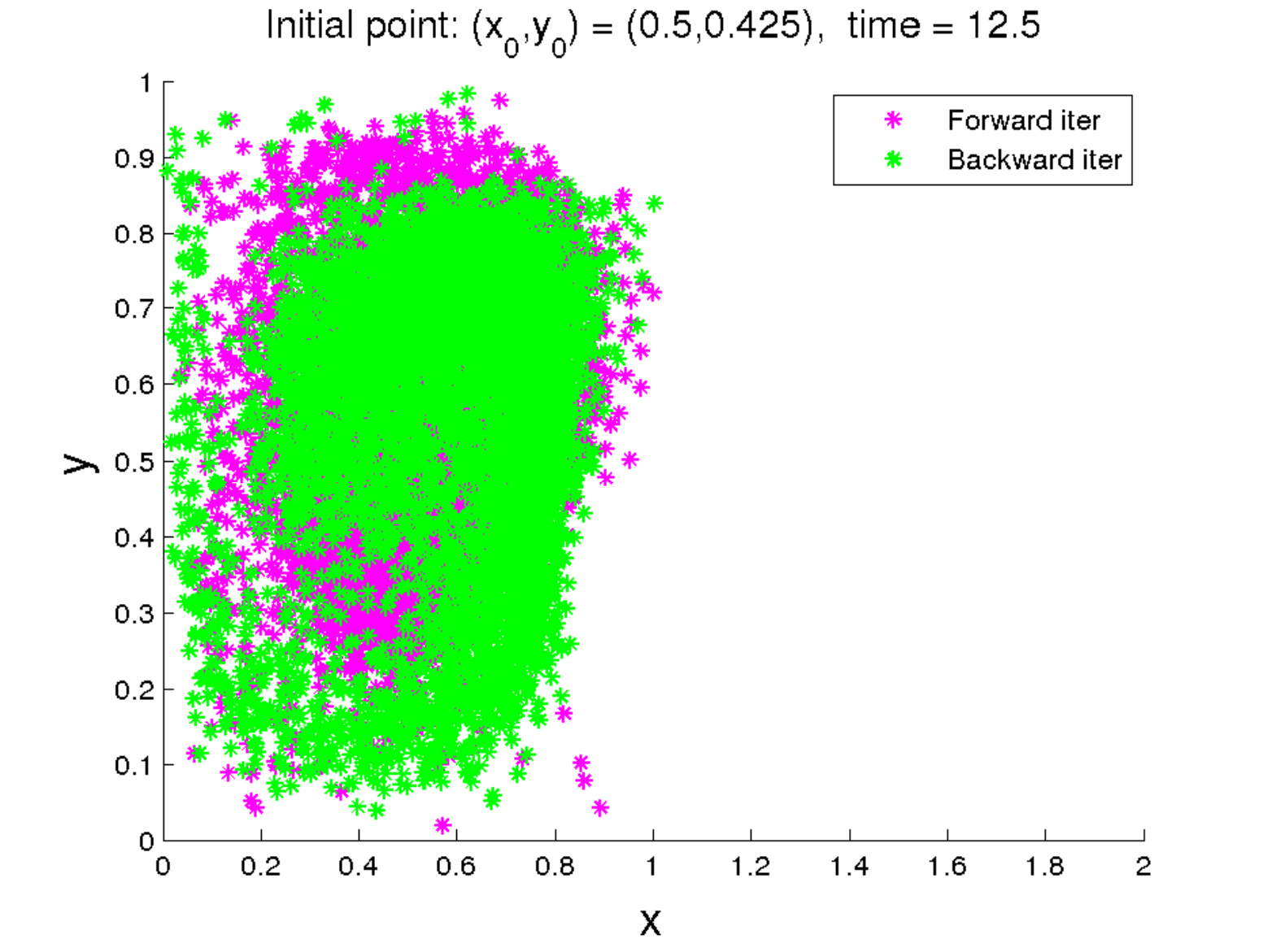}}
  \subfigure{\includegraphics[width=0.48\linewidth]{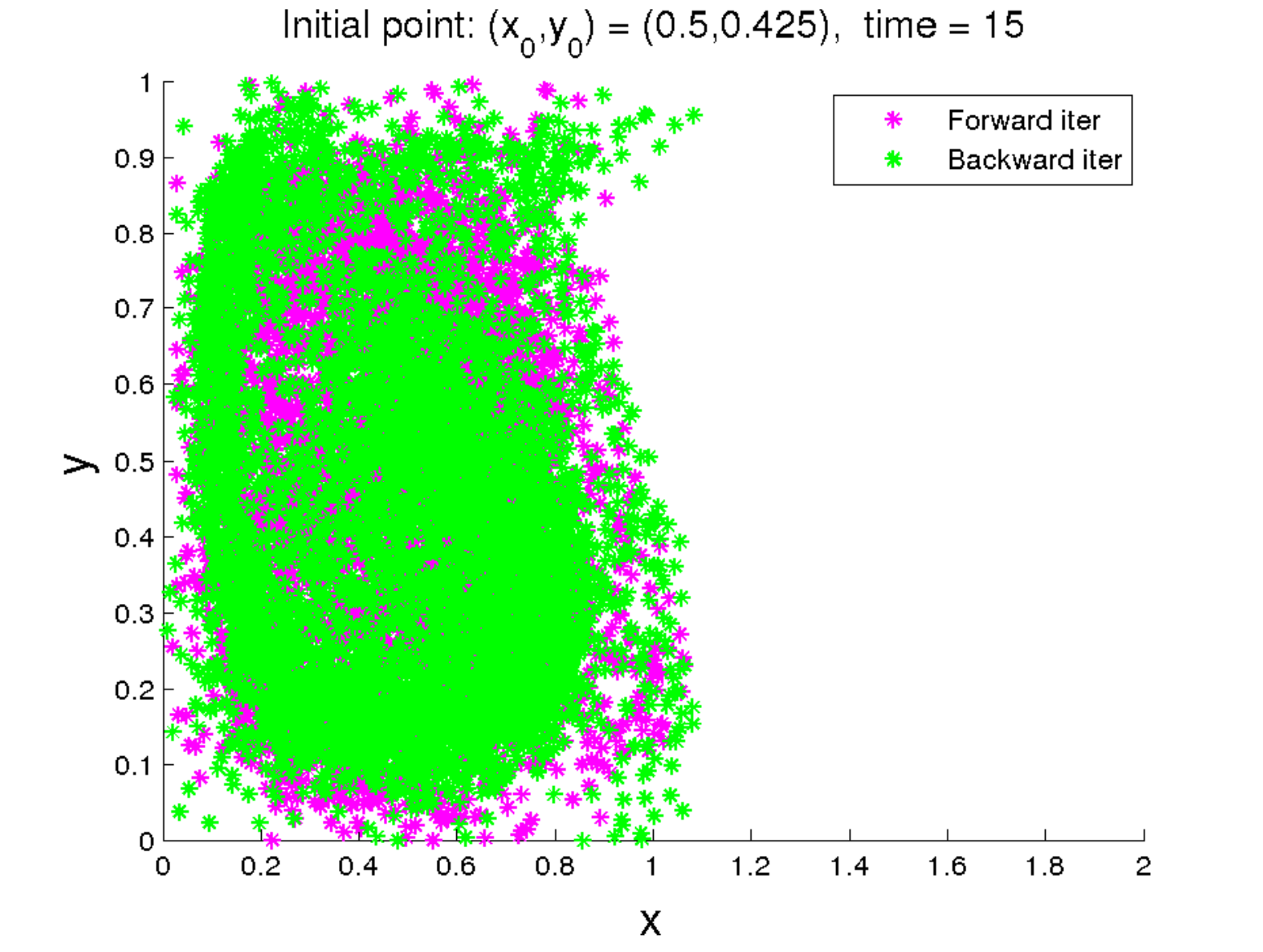}} \\
  \caption{5000 distinct simulations of trajectories starting at the point $(0.5,0.425)$ and performed 
forwards and backwards for several intervals of time.}
\label{fig:dgyre_gyre}
\end{figure*}

 
\noindent
In order to clarify the structures highlighted in Figure \ref{fig:30experiments}, we select a point in 
the middle of the gyre, the point $(0.5,0.425)$, and we evolve it forwards and backwards in time for different 
realizations of the random variable $\omega$. In Figure \ref{fig:dgyre_gyre} there are represented different snapshots 
for several units of time. As the random dynamical system evolves, the different trajectories starting at the point 
$(0.5,0.425)$ remain bounded inside the gyre. Only when the system evolves for a enough period of time, it is 
observable that hyperbolicity affects most of the trajectories.

\par
\begin{figure}[htbp!]
\centering
{\includegraphics[scale = 0.5]{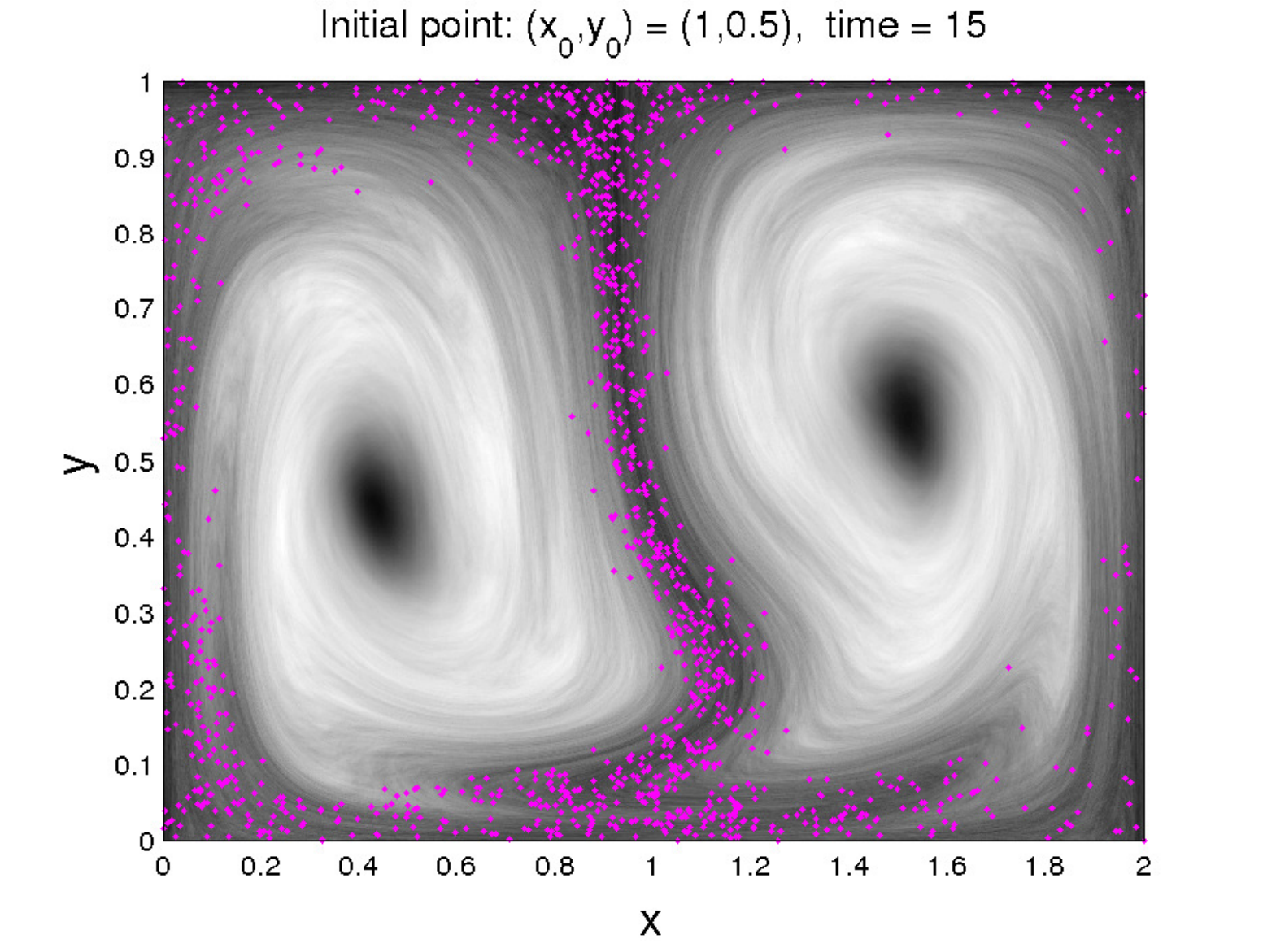}}
{\includegraphics[scale = 0.5]{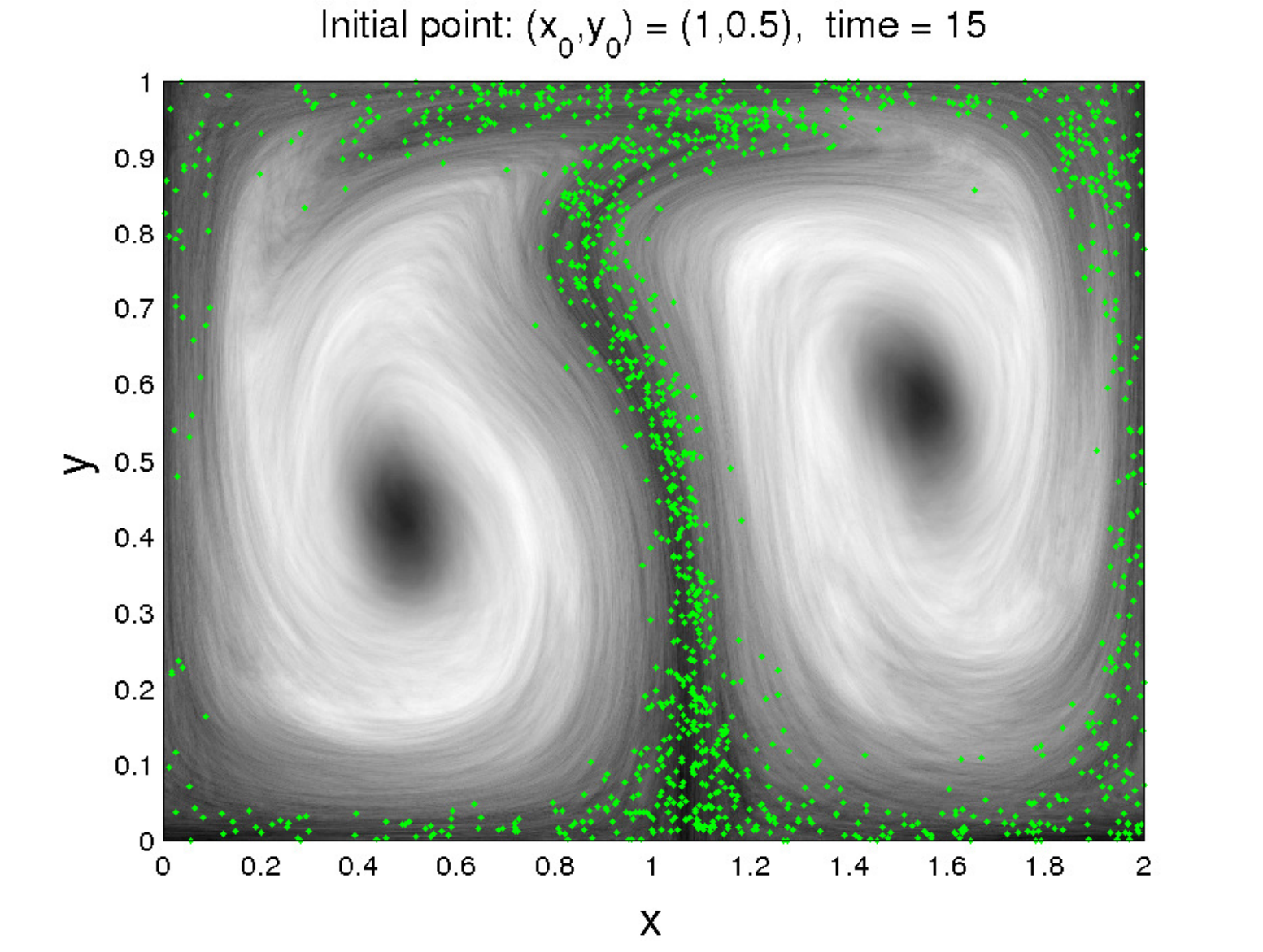}}
\caption{Representation of the ending points of 5000 distinct trajectories iterated in forward time 
from $t=0$ until $t=15$ (left hand panel, in magenta dots) and the ending points of 5000 distinct trajectories 
iterated in backward time from $t=0$ until $t=-15$ (right hand panel, in green dots). All these trajectories start at 
the same initial point $(x_{0},y_{0})=(1,0.5)$ but their respective evolutions are influenced by different Wiener 
processes. In addition, $\mathbb{E} \left[ MS_p(\cdot, \omega) \right] $ is computed for the time periods $[-\tau ,0]$ 
and $[0,\tau ]$. These expected phase spaces are depicted in the background of the left and right hand panels 
respectively with $p = 0.5$ for 30 different experiments.}
\label{fig:dgyre_inv_man}
\end{figure}

\noindent
Furthermore, in order to verify the shape that the invariant manifolds take, the point $(1,0.5)$ is evolved forward 
and backward in time for different realizations of the random variable $\omega$. In Figure \ref{fig:dgyre_inv_man} we 
overlap 5000 different trajectories starting at the same point with the expected phase space shown in the background and examine how these trajectories 
fit the sharp lines corresponding to the already computed invariant manifolds.

These pictures provide evidence for the strong correspondance between the paths followed by the particles and the structures 
depicted in bluish colours in Figure \ref{fig:30experiments}. Notice that in this last Figure \ref{fig:dgyre_inv_man}, 
the mean $\mathbb{E} \left[ MS_p(\cdot, \omega) \right] $ computed for the time period $[-\tau , \tau ]$ and represented 
in Figure \ref{fig:30experiments} has been split into two quantities. From the deterministic Lagrangian descriptors 
setting (see \cite{LBGWM15}), the trajectories iterated in forward time and starting in a neighborhood of a hyperbolic 
point approximate  the location of the unstable manifold for a large integration time $\tau$. Similarly,  the 
trajectories iterated in backward time follow the path of the stable manifold in a neighborhood of the hyperbolic point. 
In summary, these simulations support the idea of the method of stochastic Lagrangian descriptors as useful for 
depicting areas where the hyperbolic trajectories with their stable and unstable manifolds are likely located.

\section{Summary and Conclusions}
\label{sec:summ}
In this paper we have extended the method of Lagrangian descriptors to stochastic differential equations (SDEs). We have shown that the resulting method provides a way of revealing the  global phase space structure of  SDEs  that is analogous the the manner in which we understand the global phase space structure of deterministic ordinary differential equations. In particular, we show that stochastic versions of hyperbolic trajectories and their stable and unstable manifolds provide barriers to transport. We apply our method to the noisy saddle, the stochastically forced Duffing equation, and the stochastic double gyre model that is a benchmark for analyzing fluid transport.

\section*{Appendix 1. Deterministic solution mappings}
\label{appendix1}

Taking as a reference the book by Kloeden and Rasmussen (\cite{kloe11}), the solution mappings will be defined over a 
temporal domain $\mathbb{T}$, which is usually taken to be $\mathbb{R}$ (continuous time), and over a spatial domain $X$ 
which is usually taken to be $\mathbb{R}^{n}$. In both the autonomous and the nonautonomous settings, the solution mappings are denoted 
by $\phi$, and these clearly correspond to the solutions of some given deterministic differential equations.

For an autonomous differential equation $\dot{x} = f(x)$, the solution mapping $\phi$ is usually referred to as a flow, less often it is referred to as dynamical system, and it
depends only on two variables: the evolution variable $t$ and the initial condition $x_{0}$, although $\phi$ might not be defined over the entire temporal and spatial domains $\mathbb{T}$ and $X$. For any fixed initial 
condition $x_{0} \in X$, the mapping $\phi (\cdot ,x_{0})$ is a solution of the differential equation. Formally it 
satisfies the properties.

\begin{definition}
A dynamical system is a continuous function $\phi : \mathbb{T} \times X \rightarrow X$ which satisfies the following 
two properties,
\begin{enumerate}
\item[(i)] \textit{Initial value condition.} $\phi (0,x_{0})=x_{0}$ for all $x_{0} \in X$, \item[(ii)] \textit{Group 
property.} $\phi (s+t,x_{0})=\phi (s,\phi (t,x_{0}))$ for all $s,t \in \mathbb{T}$ and $x_{0} \in X$.
\end{enumerate}
\end{definition}

The solution mappings for nonautonomous differential equations of the kind $\dot{x} = f(x,t)$ are known as processes, and share 
many similarities with dynamical systems in the autonomous case. The main difference is that $\phi$ depends on an 
extra variable, the initial time $t_{0}$ when the solution $x(t)$ of the differential equation passes through $x_{0}$, 
let say $x(t_{0})=x_{0}$. This is due to the time-dependence of the vector field field $f(x,t)$.

\begin{definition}
\label{process}
A process is a continuous mapping $(t,t_{0},x_{0}) \mapsto \phi (t,t_{0},x_{0}) \in X$ for $t,t_{0} \in \mathbb{T}$ and 
$x_{0} \in X$, which satisfies the initial value and evolution properties,

\begin{enumerate}
\item[(i)] $\phi (t_{0},t_{0},x_{0})=x_{0}$ for all $t_{0} \in \mathbb{T}$ and $x_{0} \in X$, \item[(ii)] $\phi 
(t_{2},t_{0},x_{0})=\phi (t_{2},t_{1},\phi (t_{1},t_{0},x_{0}))$ for all $t_{0},t_{1},t_{2} \in \mathbb{T}$ and $x_{0} 
\in X$.
\end{enumerate}
\end{definition}

\section*{Appendix 2. RDSs for time-dependent SDEs}
\label{appendix2}

As mentioned in Section \ref{sec:pc}, just after the introduction of the formal definition of random dynamical system, this was only appropriate for SDEs with time-independent coefficients $b,\sigma$. Inspired in this Definition \ref{rds} given by J. Duan (\cite{duan15}), we develop the same notion of RDS for a more general context, where at least $b$ or $\sigma$ do depend explicitly on the time variable $t$. For this purpose it is necessary to consider a fourth argument for the solution mapping $\varphi$, which is the initial time $t_{0}$ when $\varphi$ passes through the initial condition $x_{0}$. This has been done in the same way that processes are built from flow mappings in the deterministic differential equations setting.

\begin{definition}
\label{nonautonomousRDS}
Let $\lbrace \theta_{t} \rbrace_{t \in \mathbb{R}}$ be a measure-preserving dynamical system defined over $\Omega$, and 
let $\varphi : \mathbb{R} \times \mathbb{R} \times \Omega \times \mathbb{R}^{n} \rightarrow \mathbb{R}^{n}$ be a 
measurable mapping such that $(t,t_{0},\cdot , x) \mapsto \varphi (t,t_{0},\omega ,x)$ is continuous for all $\omega \in 
\Omega$, and the family of functions $\lbrace \varphi (t,t_{0},\omega ,\cdot ): \mathbb{R}^{n} \rightarrow \mathbb{R}^{n} 
\rbrace$ has the cocycle property,
$$ \varphi (t_{0},t_{0},\omega ,x)=x \quad \text{and}$$
$$ \varphi (t_{2},t_{0},\omega ,x) = \varphi(t_{2},t_{1},\theta_{(t_{1}-t_{0})}\omega,\varphi (t_{1},t_{0},\omega ,x)) 
\quad \text{for all } t_{0},t_{1},t_{2} \in \mathbb{R}, \text{ } x \in \mathbb{R}^{n} \text{ and } \omega \in \Omega .$$
Then the mapping $\varphi$ is a random dynamical system with respect to the stochastic differential equation
$$dX_{t} = b(X_{t},t)dt + \sigma (X_{t},t)dW_{t}$$
if $\varphi (t,t_{0},\omega ,x)$ is a solution of the equation.
\end{definition}

\section*{Appendix 3. Multiplicative Ergodic  Theorem}
\label{appendix3}

The next results are taken from the book by J. Duan (2015) (\cite{duan15}). These are essential in order to understand 
the hyperbolicity of a given stationary orbit $\tilde{X}(\omega )$ with respect to a random dynamical system $\varphi$.

We begin by centering  $\varphi$ at $\tilde{X}(\omega )$, therefore achieving a new RDS,

$$\widetilde{\varphi }(t,\omega ,x) = \varphi (t,\omega ,\tilde{X}(\omega )+x) - \tilde{X}(\theta_{t}\omega ).$$

\noindent
This mapping $\widetilde{\varphi }$ retains the same dynamics as $\varphi$, with the difference that the constant value trajectory $X \equiv 0$ is a random fixed point for $\widetilde{\varphi}$. Now $\widetilde{\varphi}$ linearized with respect to the initial condition variable $x\in \mathbb{R}^{n}$ and evaluated at the initial condition $x=0$,

$$\Phi (t, \omega ) = \frac{\partial }{\partial x} \widetilde{\varphi} (t,\omega ,x=0) = \frac{\partial }{\partial x} \varphi (t,\omega ,x=\tilde{X}(\omega ))$$

\noindent
allows us to conveniently study whether  or not $\tilde{X}(\omega )$ is hyperbolic by applying the following theorem.

\begin{theorem}\textbf{Multiplicative ergodic theorem (MET).}
Let $\Phi (t, \omega )$ be a linear RDS (linear cocycle) in $\mathbb{R}^n$, for $t \in \mathbb{R}$, on a probability 
space $(\Omega , \mathcal{F}, \mathcal{P})$, over a measurable driving flow $\theta_t$. Assume that the following 
integrability conditions are satisfied:
\begin{equation}
\displaystyle{\sup_{0\leq t\leq 1}\ln^{+}\norm{\Phi (t, \omega )} \in L^1(\Omega ), \quad \sup_{0\leq t\leq 
1}\ln^{+}\norm{\Phi (-t, \omega )} \in L^1(\Omega) ,}
\label{eq:oseledets_conditions}
\end{equation}
where $\ln^{+}(z) := \max \{ \ln (z), 0\}$, denoting the nonnegative part of the natural logarithm. Then there exists 
an invariant set $\widetilde{\Omega} \in \mathcal{F}$ of full probability measure, such that for every $\omega \in 
\widetilde{\Omega}$,
\begin{enumerate}
\item[(i)] the asymptotic geometric mean $\lim_{t \to \pm \infty}[\Phi (t, \omega )^T\Phi (t, \omega )]^\frac{1}{2t} = 
\widetilde{\Phi}$ exists and it is nonnegative definite $n \times n$ matrix.
\item[(ii)] the matrix $\widetilde{\Phi}$ has distinct eigenvalues $e^{\lambda_p(\omega )} < \cdots , 
e^{\lambda_1(\omega )}, \omega \in \widetilde{\Omega}$, with corresponding eigenspaces $E_{p(\omega )}(\omega ), \cdots, 
E_1(\omega )$ of dimensions $d_i(\omega ) = \dim E_i(\omega ), i = 1, \cdots, p(\omega )$. These eigenspaces are such 
that $E_1(\omega ) \oplus E_2(\omega ) \oplus \cdots \oplus E_{p(\omega )}(\omega ) = \mathbb{R}^n$. Moreover, $p$, 
$\lambda_i$ and $d_i$ are invariant under the driving flow $\theta_t$ in the following sense:
    $$p(\theta_t \omega ) = p(\omega ), \quad \lambda_i(\theta_t \omega ) = \lambda_i(\omega ) \quad \text{ and } \quad 
d_i(\theta_t \omega ) = d_i(\omega )$$
    for $i = 1, \cdots, p, t \in \mathbb{R}$ and $\omega \in \widetilde{\Omega}$.
\item[(iii)] each $E_i(\omega )$ is invariant for the linear RDS: $\Phi (t, \omega )E_i(\omega ) = E_i(\theta_t \omega 
)$, for all $\omega \in \widetilde{\Omega}$ and all $t \in \mathbb{R}$.
\item[(iv)] $\lim_{t \to \pm \infty} 1/t \ln \norm{\Phi (t, \omega )x} = \lambda_i$ if and only if $x \in E_i(\omega ) 
- \{0\}$, for all $\omega \in \widetilde{\Omega}$ and $i = 1, \cdots, p$.
\end{enumerate}
\end{theorem}

\begin{definition}
Let $\Phi (t, \omega )$ be a linear RDS in $\mathbb{R}^n$ that satisfies the initial conditions of the MET. The vector 
subspaces based on the Oseledets spaces

$$E^{s}(\omega ) = \bigoplus_{\lambda_{i}<0}E_{i}(\omega ) \quad , \quad E^{c}(\omega ) = 
\bigoplus_{\lambda_{i}=0}E_{i}(\omega ) \quad , \quad E^{u}(\omega ) = \bigoplus_{\lambda_{i}>0}E_{i}(\omega )$$

\noindent
are called the stable, center and unstable subspaces of $\Phi (t, \omega )$, respectively. The family of its Lyapunov 
exponents and their corresponding multiplicities

$$\lbrace \lambda_{1}, \cdots , \lambda_{p} ; d_{1}, \cdots , d_{p} \rbrace$$

\noindent
is called the Lyapunov spectrum of $\Phi (t, \omega )$. Moreover if all the Lyapunov exponents are nonzero then the 
linear RDS $\Phi (t, \omega )$ is said to be hyperbolic.
\end{definition}

\vspace{0.3cm}
As mentioned before in Section \ref{sec:pc}, the random dynamical system (\ref{noisy_saddle_RDS}) associated to the noisy saddle (\ref{noisy_saddle}) satisfies the conditions of the Multiplicative ergodic theorem. This fact and the simplicity of the noisy saddle, allows us to compute explicitly the Lyapunov spectrum and the stable and unstable subspaces of the RDS $\Phi (t,\omega )$ linearized over the stationary orbit $\tilde{X}(\omega )$. 
\\
\begin{example}{(Noisy saddle point)}\\

The random dynamical system $\varphi$ corresponding to the noisy saddle (\ref{noisy_saddle}) has the expression,

$$\begin{array}{ccccccccc} \varphi : & & \mathbb{R} \times \Omega \times \mathbb{R}^{2} & & \longrightarrow & & \mathbb{R}^{2} & & \\ & & (t,\omega ,(x,y)) & & \longmapsto & & \left( \varphi_{1}(t,\omega ,x),\varphi_{2}(t,\omega ,y)\right) & = & \left( e^{t} \left( x + \int_{0}^{t}e^{-s}dW_{s}^{1}(\omega ) \right) , e^{-t} \left( y + \int_{0}^{t}e^{s}dW_{s}^{2}(\omega ) \right) \right) . \end{array} $$

\noindent
The Jacobian matrix $D_{(x,y)}\varphi$ does not depend on the third argument, let say the initial condition $(x,y)$,

$$\Phi (t) = D_{(x,y)} \varphi (t,\omega ,(x,y)) = \left( \begin{array}{cc} e^{t} & 0 \\ 0 & e^{-t} \end{array} \right) .$$

As this linear cocycle is a diagonal matrix, its norm corresponds to the largest element of its diagonal, which is also the largest eigenvalue. In this case,

$$ || \Phi (t) || = || \Phi (-t) || = \begin{cases} e^{t} \quad \text{for } t \geq 0 \\ e^{-t} \quad \text{for } t<0 \end{cases} \text{and} \quad \displaystyle{\sup_{0\leq t\leq 1}\ln^{+}\norm{\Phi (t)}} = \displaystyle{\sup_{0\leq t\leq 1}\ln^{+}\norm{\Phi (-t)}} = \displaystyle{\sup_{0\leq t\leq 1}|t|} = 1.$$

\noindent
As any constant function belongs to the space of functions $L^{1}(\Omega )$ (remember that the probability space $\Omega$ has measure equal to 1), the conditions of the MET are satisfied. Now for the asymptotic geometric mean,

$$\widetilde{\Phi} = \lim_{t \to \pm \infty}[\Phi (t )^T\Phi (t )]^\frac{1}{2t} = \left[ \left( \begin{array}{cc} e^{t} & 0 \\ 0 & e^{-t} \end{array} \right) \cdot \left( \begin{array}{cc} e^{t} & 0 \\ 0 & e^{-t} \end{array} \right) \right]^{\frac{1}{2t}} = \left( \begin{array}{cc} e^{2t} & 0 \\ 0 & e^{-2t} \end{array} \right)^{\frac{1}{2t}} = \left( \begin{array}{cc} e^{1} & 0 \\ 0 & e^{-1} \end{array} \right) ,$$

\noindent
its associated eigenvalues are $e^{1}$ and $e^{-1}$, with corresponding eigenspaces $\text{span} \lbrace \left( \begin{array}{c} 1 \\ 0 \end{array} \right) \rbrace$ and $\text{span} \lbrace \left( \begin{array}{c} 0 \\ 1 \end{array} \right) \rbrace$ respectively. These two linear subspaces are precisely the unstable $E^{u}(\omega )$ and the stable $E^{s}(\omega )$ subspaces of $\Phi (t)$, for which the Lyapunov spectrum is $\lbrace \lambda_{1}=1,\lambda_{2}=-1 \rbrace$, therefore confirming the hyperbolic nature of the stationary orbit. It results evident that for any other more complex SDE, the computation of the asymptotic geometric mean $\tilde{\Phi}$ would be more tricky; highlighting that this process is not straightforward for most  SDEs. 
\newline
\\
What remains is to establish the relationship between the linear subspaces $E^{u}(\omega )$, $E^{s}(\omega )$ and the already mentioned $\mathcal{U}(\omega )$, $\mathcal{S}(\omega )$. This is discussed in Theorem 3.1 from Mohammed and Scheutzow (\cite{moham99}) where the local unstable and the local stable manifolds of a stationary orbit $\tilde{X}(\omega )$  are related to the sets of points (in a neighborhood of $\tilde{X}(\omega )$) which are attracted by the stationary orbit in negative or positive time, respectively. 
In summary this fact presents many similarities with the definition of exponential dichotomy for deterministic dynamical systems. Indeed this exponentially attracting rate over the points of the manifolds $\mathcal{U}(\omega )$, $\mathcal{S}(\omega )$ of the noisy saddle (\ref{noisy_saddle}) is easily observed in Equation (\ref{dichotomy}),

$$(\overline{x}_{t},\overline{y}_{t}) - (\tilde{x}(\theta_{t}\omega ),\tilde{y}(\theta_{t}\omega )) = \varphi (t,\omega ,(\overline{x}_{0},\overline{y}_{0})) - \varphi (t,\omega ,(\tilde{x}(\omega ),\tilde{y}(\omega ))) = \left( \epsilon_{1}(\omega )e^{t},\epsilon_{2}(\omega )e^{-t} \right) ,$$

\noindent
regardless the initial distance $\left( \epsilon_{1}(\omega ),\epsilon_{2}(\omega ) \right)$ of $(\bar{x_{t}},\bar{y_{t}})$ to the stationary orbit $\tilde{X}(\omega )$. Moreover, as this distance can be arbitrarily large, the dynamics in a neighborhood of $\tilde{X}(\omega )$ are also achieved for the rest of trajectories within the invariant sets $\mathcal{U}(\omega )$, $\mathcal{S}(\omega )$, giving to these a global property. This argument enables us to simply refer to $\mathcal{U}(\omega )$ and $\mathcal{S}(\omega )$ as the unstable and stable manifolds of the stationary orbit of the noisy saddle equation (\ref{noisy_saddle}).
\end{example}

\vspace{0.3cm}

\section*{\bf Acknowledgments} The research of FB-I, CL and AMM is supported by the MINECO under grant MTM2014-56392-R. 
The research of SW is supported by  ONR Grant No.~N00014-01-1-0769.  We acknowledge support from MINECO: ICMAT Severo 
Ochoa project SEV-2011-0087.


\begin{thebibliography}{10}

\bibitem{arno98}
L.~Arnold.
\newblock {Random Dynamical Systems}.
\newblock {\em Springer}, 1998.

\bibitem{banisch16}
R.~Banisch and P.~Koltai.
\newblock Understanding the geometry of transport: diffusion maps for
  {L}agrangian trajectory data unravel coherent sets.
\newblock {\em http://arxiv.org/abs/1603.04709}, 2016.

\bibitem{boxl89}
P.~Boxler.
\newblock {A Stochastic Version of Center Manifold Theory.}
\newblock {\em Probability Theory and Related Fields}, 83:509--545, 1989.

\bibitem{cheng16}
Z.~Cheng, J.~Duan, and L.~Wang.
\newblock {Most probable dynamics of some nonlinear systems under noisy
  fluctuations.}
\newblock {\em Commun. Nonlinear Sci. Numer. Simulat.}, 30(1):108--114, 2016.

\bibitem{CH15}
G.~T. Craven and R.~Hernandez.
\newblock Lagrangian descriptors of thermalized transition states on
  time-varying energy surfaces.
\newblock {\em Phys. Rev. Lett.}, 115:148301, 2015.

\bibitem{CH16}
G.~T. Craven and R.~Hernandez.
\newblock Deconstructing field-induced ketene isomerization through
  {L}agrangian descriptors.
\newblock {\em Phys. Chem. Chem. Phys.}, 18:4008, 2016.

\bibitem{datta01}
S.~Datta and J.~K. Bhattacharjee.
\newblock {Effect of stochastic forcing on the Duffing oscillator.}
\newblock {\em Physics Letters A}, 283:323--326, 2001.

\bibitem{amism11}
A.~de~la C{\'a}mara, A.~M. Mancho, K.~Ide, E.~Serrano, and C.R. Mechoso.
\newblock {Routes of transport across the {A}ntarctic polar vortex in the
  southern spring}.
\newblock {\em J. Atmos. Sci.}, 69(2):753--767, 2012.

\bibitem{alvaro2}
A.~de~la C{\'a}mara, C.~R. Mechoso, E.~Serrano, and K.~Ide.
\newblock {Quasi-horizontal transport within the Antarctic polar night vortex:
  Rossby wave breaking evidence and Lagrangian structures.}
\newblock {\em J. Atmos. Sci.}, 70:2982--3001, 2013.

\bibitem{duan15}
J.~Duan.
\newblock {An Introduction to Stochastic Dynamics}.
\newblock {\em Cambridge University Press}, 2015.

\bibitem{GMWM15}
V.~J. Garc\'{i}a-Garrido, A.~M. Mancho, S.~Wiggins, and C.~Mendoza.
\newblock A dynamical systems approach to the surface search for debris
  associated with the disappearance of flight {M}{H}370.
\newblock {\em Nonlin. Proc. Geophys.}, 22:701--712, 2015.

\bibitem{GRMCW15}
V.~J. Garc\'ia-Garrido, A.~Ramos, A.~M. Mancho, J.~Coca, and S.~Wiggins.
\newblock A dynamical systems perspective for a real-time response to a marine
  oil spill.
\newblock {\em submitted to Marine Pollution Bulletin}, 2016.

\bibitem{hsieh12}
M.A. Hsieh, E.~Forgoston, T.~W. Mather, and I.B. Schwartz.
\newblock {Robotic manifold tracking of coherent structures in flows Robotics
  and Automation.}
\newblock {\em Robotics, IEEE International Conference}, 30:593--603, 2012.

\bibitem{kayo}
K.~Ide, D.~Small, and S.~Wiggins.
\newblock {Distinguished hyperbolic trajectories in time dependent fluid flows:
  analytical and computational approach for velocity fields defined as data
  sets}.
\newblock {\em Nonlin. Proc. Geophys.}, 9:237--263, 2002.

\bibitem{JH2016}
A.~Junginger and R.~Hernandez.
\newblock Lagrangian descriptors in dissipative systems.
\newblock {\em Phys. Chem. Chem. Phys.}, 2016.

\bibitem{JH16}
A.~Junginger and R.~Hernandez.
\newblock Uncovering the geometry of barrierless reactions using {L}agrangian
  descriptors.
\newblock {\em J. Phys. Chem. B}, 2016.

\bibitem{kloeden92}
P.~E. Kloeden and E.~Platen.
\newblock {\em Numerical Solution of Stochastic Differential Equations}.
\newblock Applications of Mathematics 23. Springer-Verlag Berlin Heidelberg, 1
  edition, 1992.

\bibitem{kloe11}
P.~E. Kloeden and M.~Rasmussen.
\newblock {\em Nonautonomous Dynamical Systems}, volume 176 of {\em
  Mathematical Surveys and Monographs}.
\newblock American Mathematical Society, Providence, Rhode Island, 2011.

\bibitem{carlos}
C.~Lopesino, F.~Balibrea, S.~Wiggins, and A.~M. Mancho.
\newblock {Lagrangian Descriptors for Two Dimensional, Area Preserving
  Autonomous and Nonautonomous Maps.}
\newblock {\em Communications in Nonlinear Science and Numerical Simulation},
  27(1-3):40--51, 2015.

\bibitem{LBGWM15}
C.~Lopesino, F.~Balibrea-Iniesta, V.~J. Garc\'{i}a-Garrido, S.~Wiggins, and
  A.~M. Mancho.
\newblock A theoretical framework for \text{Lagrangian} descriptors.
\newblock {\em submitted to Comm. Nonlin. Sci. Numer. Simul.}, 2016.

\bibitem{chaos}
J.~A.~J. Madrid and A.~M. Mancho.
\newblock {Distinguished trajectories in time dependent vector fields}.
\newblock {\em Chaos}, 19:013111, 2009.

\bibitem{MCWGM15}
A.~M. Mancho, J.~Curbelo, S.~Wiggins, V.~J. Garc\'{i}a-Garrido, and C.~Mendoza.
\newblock {\em Beautiful Geometries Underlying Ocean Nonlinear Processes},
  pages 80--85.
\newblock European Geophysical Union, 2015.

\bibitem{cnsns}
A.~M. Mancho, S.~Wiggins, J.~Curbelo, and C.~Mendoza.
\newblock {Lagrangian Descriptors: A Method for Revealing Phase Space
  Structures of General Time Dependent Dynamical Systems.}
\newblock {\em Communications in Nonlinear Science and Numerical Simulation},
  18:3530--3557, 2013.

\bibitem{prl}
C.~Mendoza and A.~M. Mancho.
\newblock {The hidden geometry of ocean flows}.
\newblock {\em Phys. Rev. Lett.}, 105(3):038501, 2010.

\bibitem{jfm}
C.~Mendoza and A.~M. Mancho.
\newblock {The {L}agrangian description of ocean flows: a case study of the
  {K}uroshio current}.
\newblock {\em Nonlin. Proc. Geophys.}, 19(4):449--472, 2012.

\bibitem{mmw14}
C.~Mendoza, A.~M. Mancho, and S.~Wiggins.
\newblock {Lagrangian Descriptors and the Assesment of the Predictive Capacity
  of Oceanic Data Sets.}
\newblock {\em Nonlin. Proc. Geophys.}, 21:677--689, 2014.

\bibitem{moham99}
S.-E.~A. Mohammed and M.~K.~R. Scheutzow.
\newblock {The Stable Manifold Theorem for Stochastic Differential Equations.}
\newblock {\em Ann. Probab.}, 27(2):615--652, 1999.

\bibitem{rempel}
E.~L. Rempel, A.~C.-L. Chian, A.~Brandenburg, P.~R. Munuz, and S.~C. Shadden.
\newblock Coherent structures and the saturation of a nonlinear dynamo.
\newblock {\em Journal of Fluid Mechanics}, 729:309--329, 2013.

\end{thebibliography}
\end{document}